\documentclass[11pt]{article}

\usepackage{amscd,amstex,amssymb}

\def\eqref#1{(\ref{#1})}
\newcommand{\goth}{\mathfrak}
\newcommand{\arrow}{{\:\longrightarrow\:}}
\newcommand{\Z}{{\Bbb Z}}
\newcommand{\C}{{\Bbb C}}
\newcommand{\R}{{\Bbb R}}

\newcommand{\6}{\partial}
\newcommand{\1}{\sqrt{-1}\:}
\newcommand{\inangles}[1]{{\langle #1\rangle}}

\newcommand{\restrict}[1]{{\left|_{{\phantom{|}\!\!}_{#1}}\right.}}

\renewcommand{\c}[1]{{\cal #1}}
\newcommand{\calo}{{\cal O}}

\let\oldtilde=\tilde
\renewcommand{\tilde}{\widetilde}
\renewcommand{\bar}{\overline}
\renewcommand{\phi}{\varphi}
\renewcommand{\epsilon}{\varepsilon}
\renewcommand{\geq}{\geqslant}
\renewcommand{\leq}{\leqslant}

\newcommand{\End}{\operatorname{End}}
\newcommand{\Id}{\operatorname{Id}}
\newcommand{\Tw}{\operatorname{Tw}}
\newcommand{\Vol}{\operatorname{Vol}}
\newcommand{\Hom}{\operatorname{Hom}}

\newcommand{\comment}[1]{{}}

\def\blacksquare{\hbox{\vrule width 4pt height 4pt depth 0pt}}
\def\endproof{\blacksquare}

\makeatletter

\@ifundefined{Bbb}
     {\newcommand{\Bbb}[1]{{\mathbb #1}}}%
{}

\newcommand{\ps@verbit}{%
  \renewcommand{\@oddhead}{%
          \scriptsize
          {Deformations of trianalytic subvarieties}
          \hfil\tiny {final version, Oct. 8 1996}}
  \renewcommand{\@evenhead}{\@oddhead}
  \renewcommand{\@oddfoot}{\hfil\thepage\hfil}
  \renewcommand{\@evenfoot}{\@oddfoot}}
 
\newcounter{Mycounter}[section]
\newcounter{lemma}[section]
\renewcommand{\thelemma}{{Lemma \thesection.\arabic{lemma}}}
\newcommand{\lemma}{%
     \setcounter{lemma}{\value{Mycounter}}
     \refstepcounter{lemma}
     \stepcounter{Mycounter}
     {\bf \thelemma:\ }}

\newcounter{claim}[section]
\renewcommand{\theclaim}{{Claim \thesection.\arabic{claim}}}
\newcommand{\claim}{%
     \setcounter{claim}{\value{Mycounter}}
     \refstepcounter{claim}
     \stepcounter{Mycounter}
     {\bf \theclaim:\ }}

\newcounter{sublemma}[section]
\newcounter{corollary}[section]
\renewcommand{\thecorollary}{{Corollary \thesection.\arabic{corollary}}}
\newcommand{\corollary}{%
     \setcounter{corollary}{\value{Mycounter}}
     \refstepcounter{corollary}
     \stepcounter{Mycounter}
     {\bf \thecorollary:\ }}

\newcounter{theorem}[section]
\renewcommand{\thetheorem}{{Theorem \thesection.\arabic{theorem}}}
\newcommand{\theorem}{%
     \setcounter{theorem}{\value{Mycounter}}
     \refstepcounter{theorem}
     \stepcounter{Mycounter}
     {\bf \thetheorem:\ }}

\newcounter{proposition}[section]
\renewcommand{\theproposition}
       {{Proposition \thesection.\arabic{proposition}}}
\newcommand{\proposition}{%
     \setcounter{proposition}{\value{Mycounter}}
     \refstepcounter{proposition}
     \stepcounter{Mycounter}
     {\bf \theproposition:\ }}

\newcounter{definition}[section]
\renewcommand{\thedefinition}
       {{Definition \thesection.\arabic{definition}}}
\newcommand{\definition}{%
     \setcounter{definition}{\value{Mycounter}}
     \refstepcounter{definition}
     \stepcounter{Mycounter}
     {\bf \thedefinition:\ }}

\newcounter{example}[section]
\renewcommand{\theexample}{{Example \thesection.\arabic{example}}}
\newcommand{\example}{%
     \setcounter{example}{\value{Mycounter}}
     \refstepcounter{example}
     \stepcounter{Mycounter}
     {\bf \theexample:\ }}

\newcounter{remark}[section]
\renewcommand{\theremark}{{Remark \thesection.\arabic{remark}}}
\newcommand{\remark}{%
     \setcounter{remark}{\value{Mycounter}}
     \refstepcounter{remark}
     \stepcounter{Mycounter}
     {\bf \theremark:\ }}

\newcounter{problem}[section]
\newcounter{question}[section]
\@addtoreset{equation}{section}
\@addtoreset{footnote}{section}
\makeatother

\begin{document}

\begin{center}
{\Large\bf
	Deformations of trianalytic subvarieties of \\[2mm] hyperk\"ahler 
manifolds.}\\[4mm]
Misha Verbitsky,\footnote{Supported by the NSF grant 9304580}\\[4mm]
{\tt verbit@@thelema.dnttm.rssi.ru, verbit@@math.ias.edu}
\end{center}

\hfill

{\small 
\hspace{0.2\linewidth}
\begin{minipage}[t]{0.7\linewidth}
Let $M$ be a compact complex manifold equipped 
with a hyperk\"ahler metric, and $X$ be a closed 
complex analytic subvariety of $M$.
In alg-geom 9403006, we proved that $X$ is {\bf trianalytic}
(i. e., complex analytic with respect to all complex structures induced
by the hyperk\"ahler structure), provided that $M$ is generic in 
its deformation class. Here we study the complex analytic deformations
of trianalytic subvarieties. We prove that all deformations of $X$ are
trianalytic and 
naturally isomorphic to $X$ as complex analytic varieties. We show
that this isomorphism is compatible with the metric induced from $M$. 
Also, we prove that the Douady space of complex analytic
deformations of $X$ in $M$ is equipped with a natural 
hyperk\"ahler structure.
\end{minipage}
}

\tableofcontents

\hfill


\section{Introduction.}
\label{_Intro_Section_}


\subsection{An overview}

This is the third article studying closed complex analytic subvarieties
of compact holomorphically symplectic\footnote{Holomorphically
symplectic means ``equipped with a holomorphic symplectic form''.
See \ref{_holomorphi_symple_Definition_} for details.}
 K\"ahler manifolds. In the first article in series 
(\cite{Verbitsky:Symplectic_I_}), we proved that, when a holomorphically
symplectic manifold $M$ is generic in its deformation class, all
subvarieties $X\subset M$ are also holomorphically symplectic, 
i. e. restriction of holomorphic symplectic form to nonsingular
strata of $X$ is non-degenerate. In the second article
\cite{Verbitsky:Symplectic_II_}, we obtained a more precise result
about the structure of such $X$, which is related to the hyperk\"ahler
structure of $M$.

By Yau's proof of Calabi conjecture, $M$ admits a natural hyperk\"ahler
metric (\ref{_symplectic_=>_hyperkahler_Proposition_};
for a definition of hyperk\"ahler manifold,
see \ref{_hyperkahler_manifold_Definition_}). 
A hyperk\"ahler structure (which is
essentially a quaternion action in the tangent bundle to $M$) gives
a rise to a whole family of complex structure on $M$, parametrized
by $\C P^1$. These complex structures are called {\bf complex
structures induced by the hyperk\"ahler structure}
(\ref{_indu_comple_str_Definition_}).
 Denote the set of all induced complex structures
by $\c R_M$. A closed subset $X\subset M$ is called 
{\bf trianalytic} if $X$ is complex analytic with respect
to all induced complex structures $L \in \c R_M$
(\ref{_trianalytic_Definition_}). For an induced
complex structure $L$, we denote by $(M, L)$ the $M$ considered
as a complex manifold, with complex structure $L$. In
\cite{Verbitsky:Symplectic_II_}, we proved that for all $L\in \c R_M$,
with exception of may be a countable set, all complex analytic subsets
of $(M, L)$ are trianalytic (\ref{_hyperkae_embeddings_Corollary_}).

Unlike the second article \cite{Verbitsky:Symplectic_II_},
which supersedes results of the first \cite{Verbitsky:Symplectic_I_},
the present one (the third) elaborates on the results of the second.
We start where \cite{Verbitsky:Symplectic_II_} left. Consider
a complex subvariety $X \subset (M, L)$ which happens to be
trianalytic. Such subvarieties are
called {\bf complex analytic subvarieties of trianalytic type}. We 
describe the deformations of $X$ in $(M, L)$ and 
the Douady space of such deformations.\footnote{The Douady space
\cite{_Douady_} is the same as Chow scheme, but in complex analytic
situation.}

In \cite{Verbitsky:Symplectic_II_}, we gave a simple cohomological
criterion of trianaliticity, for arbitrary complex analytic
subvariety $X \subset (M, L)$ 
(see \ref{_G_M_invariant_implies_trianalytic_Theorem_} 
of this article). \ref{_G_M_invariant_implies_trianalytic_Theorem_}  
immediately implies that a complex analytic deformation of 
a trianalytic subvariety is again trianalytic.
The main result of this article is the following theorem.

\hfill

\theorem \label{_iso_intro:Theorem_}
Let $M$ be a hyperk\"ahler manifold, $L$ an induced complex structure,
and $X, X'\subset (M, L)$ be closed complex analytic subvarieties in
the same deformation class. Assume that $X$ is
trianalytic.\footnote{Then the subvariety $X'$ is also trianalytic,
as implied by \ref{_G_M_invariant_implies_trianalytic_Theorem_}.} 
Consider $X, X'$ as K\"ahler subvarieties of $(M, L)$,
with K\"ahler metric induced from $M$. Then 

\begin{description}
\item[(i)] $X, X'$ are naturally
isomorphic as complex varieties, and this isomorphism is compatible
with the K\"ahler metric.

\item[(ii)] The embedding $X\stackrel i \hookrightarrow M$ 
is compatible with metrics given by geodesics. 

\item[(iii)] Also, every map 
$X\stackrel \phi\arrow (M, L)$
is complex analytic, provided that $\phi$ is a 
deformation of $i$ in the space of isometries
$X\arrow M$.
\end{description}

\hfill

{\bf Caution} The varieties $X, X'$ need not to be non-singular.

\hfill

\ref{_iso_intro:Theorem_} (i) follows from
\ref{_triana_subse_isome_Theorem_} and
\ref{_triana_subse_comple_ana_Theorem_}.
\ref{_iso_intro:Theorem_} 
(ii) follows from \ref{_hype_embe_comple_geode:Corollary_},
and (iii) from \ref{_isome_embe_Proposition_}.

\hfill

As a corollary, we obtain the following interesting result
(\ref{_Doua_hyperka_Theorem_}).

\hfill

\theorem \label{_defo_intro:Theorem_} 
Let $M$ be a compact holomorphically symplectic K\"ahler 
manifold, and $X\subset M$ a closed complex analytic subvariety.
Consider the Douady space $D_M(X)$ of deformations of $X$ inside
$M$. Then $D_M(X)$ is compact and is equipped with a natural
hyperk\"ahler structure, in the sense of 
\cite{_Verbitsky:Hyperholo_bundles_},\footnote{If $D_M(X)$
is non-singular, this means exactly that $D_M(X)$ is hyperk\"ahler.
In singular case, there is no satisfactory definition of a hyperk\"ahler
structure. The definition of \cite{_Verbitsky:Hyperholo_bundles_}
is a palliative, which is probably much stronger that the correct
definition. See \cite{_Simpson:hyperka-defi_}, \cite{_Deligne:defi_}
for alternative definitions.} (see \ref{_singu_hype_Definition_}).

\hfill

\remark This gives a new set of examples of hyperk\"ahler
varieties, in addition to those produced by 
\cite{_Verbitsky:Hyperholo_bundles_}. The varieties obtained
through \cite{_Verbitsky:Hyperholo_bundles_} are usually not
compact; on contrary, our new examples are compact.

\hfill

The proof of \ref{_iso_intro:Theorem_} is based on the
following argument. First, we show that every trianalytic
submanifold $X\subset M$ (not necessarily closed) is 
{\bf completely geodesic} in $M$ (\ref{_hype_embe_comple_geode:Corollary_}).
\footnote{{\bf Completely geodesic} means 
that geodesics in $X$ are also geodesics in $M$.}
A simple geometric argument shows that a family $\c X$ of completely
geodesic submanifolds admits a natural {\bf connection}
(\ref{_conne_Definition_},  
\ref{_conne_in_fam_of_comple_geode_Proposition_}).
This connection is compatible with the holomorphic
structure if the family $\c X$ is holomorphic
(\ref{_conne_in_fam_of_comple_geode_Proposition_}).

\hfill

For $M$ a compact hyperk\"ahler variety,
$I$ an induced complex structure and $X$ a trianalytic subvariety of $M$,
we study the Douady space $D(X,I)$ which classifies complex analytic
deformations of $X$ in $(M,I)$, where $(M, I)$ is $M$
considered as a complex manifold. We prove that
the real analytic subvariety $D(X)$ underlying
$D(X,I)$ does not depend from the choice of induced
complex structure. This endows $D(X)$ with a 
2-dimensional sphere of complex structures, which
induce quaternionic action in the Zariski tangent space. 
A real analytic variety with such a system of complex
structures is called {\bf hypercomplex} 
(\ref{_hypercomplex_Definition_}).

We return to the study of the families of trianalytic subvarieties.
As we have shown, the base of the universal family is 
hypercomplex. The fibers are trianalytic, and therefore,
hypercomplex as well. The
natural connection in such a family is compatible with
the quaternionic action, because this connection is 
{\bf holomorphic} with respect to each of 
induced complex structures. Thus, a curvature of this
connection is ${\Bbb H}^*$-invariant (\ref{_curva_SU(2)_inva_Lemma_}). 
On the other hand, the curvature lies in the tensor product
of three representations of ${\Bbb H}^*$ of weight one.
A trivial linear algebra argument shows that such a
tensor product does not contain non-trivial ${\Bbb H}^*$-invariant
vectors. This proves that the natural connection
in the family of trianalytic subvarieties is {\bf flat}
(\ref{_conne_in_triana_flat_Theorem_}).

Let $\c X\stackrel \pi \arrow S$ be a family of closed trianalytic
subvarieties in a compact hyperk\"ahler manifold $M$. Let
$X_s^{ns}$ be the set of non-singular points in
$X_s:= \pi^{-1}(s)$, $s\in S$. The submanifold $X_s^{ns}$
is completely geodesic. Thus, the family $\c X$ is equipped
with a natural connection. For every two points $s_1, s_2 \in S$,
we show that this connection might be integrated to an
isomorphism of hyperk\"ahler manifolds
$\Psi^{s_2}_{s_1}:\; X^{ns}_{s_1} \arrow X^{ns}_{s_2}$
(\ref{_triana_subse_isome_Theorem_}). Since $X_s^{ns}$ is
completely geodesic in $M$, its completion as a metric space
is naturally isomorphic to a closure $X_s$ of $X_s^{ns}$ in $M$.
Since the map $\Psi^{s_2}_{s_1}:\; X^{ns}_{s_1} \arrow X^{ns}_{s_2}$
is an isometry, it is naturally extended to an isomorphism
of metric spaces $\bar \Psi^{s_2}_{s_1}:\; X^{ns}_{s_1} \arrow
X^{ns}_{s_2}$ (\ref{_triana_subse_isome_Theorem_}). From
construction it follows that $\Psi^{s_2}_{s_1}$ is compatible
with each of induced holomorphic structures. We extend
this assertion to $\bar \Psi^{s_2}_{s_1}$. This is done
in two steps. A general type argument shows that
a homeomorphism of complex varieties is bimeromorphic
if it is holomorphic in a dense open set. The leap from
bimeromorphic to holomorphic is done via a convoluted
algebro-geometric argument involving normalization
and finite unramified maps.

We finish this paper with the definition of singular hyperk\"ahler
varieties. We cite the previously known examples of singular hyperk\"ahler
varieties and then show that constructed above maps 
$\bar \Psi^{s_2}_{s_1}:\; X^{ns}_{s_1} \arrow X^{ns}_{s_2}$
are isomorphisms of hyperk\"ahler varieties. The Douady space
of trianalytic subvarieties of a compact hyperk\"ahler manifold
is proven to be hyperk\"ahler.

\subsection{Contents}

\begin{itemize}

\item Section \ref{_Intro_Section_} is an introduction,
independent from the rest of this paper.

\item Section \ref{_basics_Section_} is a compendium of results
pertaining to hyperk\"ahler geometry and Yang--Mills theory.
We define hyperk\"ahler manifolds, trianalytic subvarieties, 
hyperholomorphic bundles, and cite their most basic properties.
A reader acquainted with a literature may skip this section.

\item Section \ref{_nonsingu_preli_Section_} illustrates our results
with a simplicistic example of deformations of a smooth trianalytic
subvariety. This section is also independent from the rest,
and its results are further superseded by
more general statements. It is safe to skip this section too.

\item Section \ref{_comple_geode_defo_Section_} contains a study
of completely geodesic submanifolds and their deformations.
Using an easy geometric argument,
we show that a family of completely geodesic submanifolds
is equipped with a natural connection
(\ref{_conne_in_fam_of_comple_geode_Proposition_}). We study this
connection and show that, under certain additional assumptions,
this connection induces isometry (metric equivalence)
of the fibers of this family. This section does not use
results of hyperk\"ahler geometry.

\item In Section \ref{_comple_geode_hyperho_Section_}, 
we prove that hyperk\"ahler
submanifolds are always completely geodesic. We use arguments from 
Yang--Mills theory and the theory of hyperholomorphic
bundles (\cite{_Verbitsky:Hyperholo_bundles_}).

\item Section \ref{_triholo_Section_} shows that the deformational results
of Section \ref{_comple_geode_defo_Section_} can be applied to the
deformations of hyperk\"ahler submanifolds, not necessary closed. 
Again, we use arguments from the theory of hyperholomorphic
bundles. In Appendix to this section, we prove equivalence of real
analytic structures induced on a hyperk\"ahler manifold by different
induced complex structures. The proof is based on twistor
geometry (\cite{_HKLR_}, \cite{_NHYM_}).

\item Section \ref{_Douady_Section_} studies the Douady deformation
space $D(X)$ of trianalytic subvariety $X$ of a compact hyperk\"ahler 
manifold $M$.
A general argument of K\"ahler geometry shows that the Douady
space is compact (\cite{_Lieberman_}). Using this, we prove that
the underlying real analytic variety does not depend on the
complex structure. As an application, we obtain that
the Douady space is {\bf hypercomplex} (\ref{_hypercomplex_Definition_}),
i. e., has an integrable quaternionic action in the tangent space.
In the Appendix to this section, we give an independent
proof of the compactness of the Douady space. We describe
the Douady space in terms of the space of isometric embeddings
from $X$ to $M$. The Appendix is based on Wirtinger's
inequality (\ref{_Wirti_for_Kahle_Theorem_}). 
The body of this section depends only on the result
of the Appendix to Section \ref{_triholo_Section_} 
(the equivalence of induced real analytic structures
on a given hyperk\"ahler manifold).

\item In Section \ref{_Conne_in_fami_Section_},
we define a curvature of a connection in a family of
manifolds. We consider the natural connection on a family
of trianalytic subvarieties of a compact hyperk\"ahler manifolds.
Using the results of Section \ref{_Douady_Section_}
(hypercomplex structure on the Douady space), we show that this connection
is {\bf flat} i. e. its curvature is zero. This result is
not used anywhere outside of Section \ref{_singu_hype_Section_}
(and even there, we don't really need it). However,
the flatness of the connection shows that the
natural isomorphism $\Psi^{s_1}_{s_2}$ of the fibers
of a family $\c X \arrow S$ of trianalytic subvarieties
is independent from the choice of a path $\gamma:\; [0,1] \arrow S$
in its homotopy class.

\item Section \ref{_isome=>holo_Section_} deals with the map 
$\Psi^{s_1}_{s_2}:\; X_{s_1} \arrow X_{s_2}$ 
of fibers of a family of trianalytic subvarieties of $M$,
obtained by integrating the natural connection. This map is by
construction a homeomorphism and is holomorphic outside of singularities,
with respect to each of induced complex structures. Using general argument,
we show that $\Psi^{s_1}_{s_2}$ is bimeromorphic and induces an isomorphism
of normalizations of $X_{s_1}$, $X_{s_2}$. Then, we apply the knowledge
that, outside of singularities,
$X_{s_i}$ are completely geodesic in $M$. This allows us to 
show that $\Psi^{s_1}_{s_2}$ induces an isomorphism on the
Zariski tangent spaces. An algebro-geometric argument shows
that these properties are sufficient to prove that
$\Psi^{s_1}_{s_2}:\; X_{s_1} \arrow X_{s_2}$ is a complex analytic
isomorphism, with respect to each of induced complex structures.

In Section \ref{_isome=>holo_Section_}, we use the results of 
Section \ref{_comple_geode_defo_Section_} (structure of deformations
of completely geodesic submanifolds) and Section 
\ref{_triholo_Section_} (that these results may be applied to
the deformations of trianalytic submanifolds). 

\item In Section \ref{_singu_hype_Section_}, we give a definition
of singular hyperk\"ahler varieties. After a short discussion,
we give examples of hyperk\"ahler varieties, stemming from the
theory of hyperholomorphic bundles (\cite{_Verbitsky:Hyperholo_bundles_}).
All trianalytic subvarieties of hyperk\"ahler manifolds are
hyperk\"ahler varieties, which is clear from the definition.
Also, the map $\Psi^{s_1}_{s_2}:\; X_{s_1} \arrow X_{s_2}$
constructed in Section \ref{_isome=>holo_Section_}, is an
isomorphism of hyperk\"ahler varieties. We show that
the Douady space of deformations of a trianalytic subvariety
is also a hyperk\"ahler variety.

\end{itemize}


\section{Basic definitions and results.}
\label{_basics_Section_}


This section used mainly for reference, and contains a 
compilation of results and definitions from literature. 
An impatient reader is advised to skip it and proceed to the next 
section.

\subsection{Hyperk\"ahler manifolds}
\label{_basics_hyperka_Section_}


This subsection contains a compression of 
the basic and most known results 
and definitions from hyperk\"ahler geometry, found, for instance, in
\cite{_Besse:Einst_Manifo_} or in \cite{_Beauville_}.

\hfill

\definition \label{_hyperkahler_manifold_Definition_} 
(\cite{_Besse:Einst_Manifo_}) A {\bf hyperk\"ahler manifold} is a
Riemannian manifold $M$ endowed with three complex structures $I$, $J$
and $K$, such that the following holds.
 
\begin{description}
\item[(i)]  the metric on $M$ is K\"ahler with respect to these complex 
structures and
 
\item[(ii)] $I$, $J$ and $K$, considered as  endomorphisms
of a real tangent bundle, satisfy the relation 
$I\circ J=-J\circ I = K$.
\end{description}

\hfill 

The notion of a hyperk\"ahler manifold was 
introduced by E. Calabi (\cite{_Calabi_}).

\hfill
 
Clearly, hyperk\"ahler manifold has the natural action of
quaternion algebra ${\Bbb H}$ in its real tangent bundle $TM$. 
Therefore its complex dimension is even.
For each quaternion $L\in \Bbb H$, $L^2=-1$,
the corresponding automorphism of $TM$ is an almost complex
structure. It is easy to check that this almost 
complex structure is integrable (\cite{_Besse:Einst_Manifo_}).

\hfill

\definition \label{_indu_comple_str_Definition_} 
Let $M$ be a hyperk\"ahler manifold, $L$ a quaternion satisfying
$L^2=-1$. The corresponding complex structure on $M$ is called
{\bf an induced complex structure}. The $M$ considered as a complex
manifold is denoted by $(M, L)$.

\hfill

\definition \label{_holomorphi_symple_Definition_} 
Let $M$ be a complex manifold and $\Omega$ a closed 
holomorphic 2-form over $M$ such that 
$\Omega^n=\Omega\wedge\Omega\wedge...$, is
a nowhere degenerate section of a canonical class of $M$
($2n=dim_\C(M)$).
Then $M$ is called {\bf holomorphically symplectic}.

\hfill

Let $M$ be a hyperk\"ahler manifold; denote the
Riemannian form on $M$ by $<\cdot,\cdot>$.
Let the form $\omega_I := <I(\cdot),\cdot>$ be the usual K\"ahler
form  which is closed and parallel
(with respect to the Levi-Civitta connection). Analogously defined 
forms $\omega_J$ and $\omega_K$ are
also closed and parallel. 
 
A simple linear algebraic
consideration (\cite{_Besse:Einst_Manifo_}) shows that the form
$\Omega:=\omega_J+\sqrt{-1}\omega_K$ is of
type $(2,0)$ and, being closed, this form is also holomorphic.
Also, the form $\Omega$ is nowhere degenerate, as another linear 
algebraic argument shows.
It is called {\bf the canonical holomorphic symplectic form
of a manifold M}. Thus, for each hyperk\"ahler manifold $M$,
and induced complex structure $L$, the underlying complex manifold
$(M,L)$ is holomorphically symplectic. The converse assertion
is also true:

\hfill

\proposition \label{_symplectic_=>_hyperkahler_Proposition_}
(\cite{_Beauville_}, \cite{_Besse:Einst_Manifo_})
Let $M$ be a compact holomorphically
symplectic K\"ahler manifold with the holomorphic symplectic form
$\Omega$, a K\"ahler class 
$[\omega]\in H^{1,1}(M)$ and a complex structure $I$. 
There is a unique hyperk\"ahler structure $(I,J,K,(\cdot,\cdot))$
over $M$ such that the cohomology class of the symplectic form
$\omega_I=(\cdot,I\cdot)$ is equal to $[\omega]$ and the
canonical symplectic form $\omega_J+\1\omega_K$ is
equal to $\Omega$.

\hfill

\ref{_symplectic_=>_hyperkahler_Proposition_} immediately
follows from the conjecture of Calabi, pro\-ven by
Yau (\cite{_Yau:Calabi-Yau_}). 
\endproof


\subsection{Hyperholomorphic bundles}
\label{_hyperholo_Subsection_}


This subsection contains several versions of a
definition of hyperholomorphic connection in a complex
vector bundle over a hyperk\"ahler manifold.
We follow \cite{_Verbitsky:Hyperholo_bundles_}.

\hfill

Let $M$ be a hyperk\"ahler manifold. We identify the group $SU(2)$
with the group of unitary quaternions. This gives a canonical 
action of $SU(2)$ on the tangent bundle, and all its tensor
powers. In particular, we obtain a natural action of $SU(2)$
on the bundle of differential forms. Of a special interest
to us are those forms which are $SU(2)$-invariant.

\hfill

\lemma \label{_SU(2)_inva_type_p,p_Lemma_} 
Let $\omega$ be a differential form over
a hyperk\"ahler manifold $M$. The form $\omega$ is $SU(2)$-invariant
if and only if it is of Hodge type $(p,p)$ with respect to all 
induced complex structures on $M$.

{\bf Proof:} This is \cite{_Verbitsky:Hyperholo_bundles_}, 
Proposition 1.2. \endproof

\hfill

Further in this article, we use the following statement.

\lemma \label{_SU(2)_commu_Laplace_Lemma_}
The action of $SU(2)$ on differential forms commutes
with the Laplacian.

{\bf Proof:} This is Proposition 1.1
of \cite{_Verbitsky:Hyperholo_bundles_}. \endproof

Thus, for compact $M$, we may speak of the natural action of
$SU(2)$ in cohomology.

\hfill

 Let $B$ be a holomorphic vector bundle over a complex
manifold $M$, $\nabla$ a  connection 
in $B$ and $\Theta\in\Lambda^2\otimes End(B)$ be its curvature. 
This connection
is called {\bf compatible with a holomorphic structure} if
$\nabla_X(\zeta)=0$ for any holomorphic section $\zeta$ and
any antiholomorphic tangent vector $X$. If there exist
a holomorphic structure compatible with the given
Hermitian connection then this connection is called
{\bf integrable}.
 
\hfill
 
One can define a {\bf Hodge decomposition} in the space of differential
forms with coefficients in any complex bundle, in particular,
$End(B)$.

\hfill

\theorem \label{_Newle_Nie_for_bu_Theorem_}
Let $\nabla$ be a Hermitian connection in a complex vector
bundle $B$ over a complex manifold. Then $\nabla$ is integrable
if and only if $\Theta\in\Lambda^{1,1}(M, \End(B))$, where
$\Lambda^{1,1}(M, \End(B))$ denotes the forms of Hodge
type (1,1). Also,
the holomorphic structure compatible with $\nabla$ is unique.

{\bf Proof:} This is Proposition 4.17 of \cite{_Kobayashi_}, 
Chapter I.
$\blacksquare$

\hfill

\definition \label{_hyperho_conne_Definition_}
Let $B$ be a Hermitian vector bundle with
a connection $\nabla$ over a hyperk\"ahler manifold
$M$. Then $\nabla$ is called {\bf hyperholomorphic} if 
$\nabla$ is
integrable with respect to each of the complex structures induced
by the hyperk\"ahler structure. 
 
As follows from 
\ref{_Newle_Nie_for_bu_Theorem_}, $\nabla$ is hyperholomorphic
if and only if its curvature $\Theta$ is of Hodge type (1,1) with
respect to any of complex structures induced by a hyperk\"ahler 
structure.

As follows from \ref{_SU(2)_inva_type_p,p_Lemma_}, 
$\nabla$ is hyperholomorphic
if and only if $\Theta$ is a $SU(2)$-invariant differential form.

\hfill

\example \label{_tangent_hyperholo_Example_} 
(Examples of hyperholomorphic bundles)

\begin{description}

\item[(i)]
Let $M$ be a hyperk\"ahler manifold, $TM$ its tangent bundle
equipped with Levi--Civita connection $\nabla$. Then $\nabla$
is integrable with respect to each induced complex structure,
and hence, Yang--Mills.

\item[(ii)] For $B$ a hyperholomorphic bundle, all its tensor powers
are also hyperholomorphic.

\item[(iii)] Thus, the bundles of differential forms on a hyperk\"ahler
manifold are also hyperholomorphic.

\end{description}


\subsection{Stable bundles and Yang--Mills connections.}


This subsection is a compendium of the most
basic results and definitions from the Yang--Mills theory
over K\"ahler manifolds, concluding in the fundamental
theorem of Uhlenbeck--Yau \cite{_Uhle_Yau_}.

\hfill

\definition\label{_degree,slope_destabilising_Definition_} 
Let $F$ be a coherent sheaf over
an $n$-dimensional compact K\"ahler manifold $M$. We define
$\deg(F)$ as
 
\[ 
   \deg(F)=\int_M\frac{ c_1(F)\wedge\omega^{n-1}}{vol(M)}
\] 
and $\text{slope}(F)$ as
\[ 
   \text{slope}(F)=\frac{1}{\text{rank}(F)}\cdot \deg(F). 
\]
The number $\text{slope}(F)$ depends only on a
cohomology class of $c_1(F)$. 

Let $F$ be a coherent sheaf on $M$
and $F'\subset F$ its proper subsheaf. Then $F'$ is 
called {\bf destabilizing subsheaf} 
if $\text{slope}(F') \geq \text{slope}(F)$

A holomorphic vector bundle $B$ is called {\bf stable}
if it has no destabilizing subsheaves.
 
\hfill

Later on, we usually consider the bundles $B$ with $deg(B)=0$.

\hfill

Let $M$ be a K\"ahler manifold with a K\"ahler form $\omega$.
For differential forms with coefficients in any vector bundle
there is a Hodge operator $L: \eta\arrow\omega\wedge\eta$.
There is also a fiberwise-adjoint Hodge operator $\Lambda$
(see \cite{_Griffi_Harri_}).
 
\hfill

\definition \label{Yang-Mills_Definition_}
Let $B$ be a holomorphic bundle over a K\"ahler manifold $M$
with a holomorphic Hermitian connection $\nabla$ and a 
curvature $\Theta\in\Lambda^{1,1}\otimes End(B)$.
The Hermitian metric on $B$ and the connection $\nabla$
defined by this metric are called {\bf Yang-Mills} if 

\[
   \Lambda(\Theta)=constant\cdot \Id\restrict{B},
\]
where $\Lambda$ is a Hodge operator and $\Id\restrict{B}$ is 
the identity endomorphism which is a section of $End(B)$.

Further on, we consider only these Yang--Mills connections
for which this constant is zero.

\hfill

A holomorphic bundle is called  {\bf indecomposable} 
if it cannot be decomposed onto a direct sum
of two or more holomorphic bundles.

\hfill

The following fundamental 
theorem provides examples of Yang-\--Mills \linebreak bundles.

\theorem \label{_UY_Theorem_} 
(Uhlenbeck-Yau)
Let B be an indecomposable
holomorphic bundle over a compact K\"ahler manifold. Then $B$ admits
a Hermitian Yang-Mills connection if and only if it is stable, and
this connection is unique.
 
{\bf Proof:} \cite{_Uhle_Yau_}. \endproof

\hfill

\proposition \label{_hyperholo_Yang--Mills_Proposition_}
Let $M$ be a hyperk\"ahler manifold, $L$
an induced complex structure and $B$ be a complex vector
bundle over $(M,L)$. 
Then every hyperholomorphic connection $\nabla$ in $B$
is Yang-Mills and satisfies $\Lambda(\Theta)=0$
where $\Theta$ is a curvature of $\nabla$.
 
\hfill

{\bf Proof:} We use the definition of a hyperholomorphic 
connection as one with $SU(2)$-invariant curvature. 
Then \ref{_hyperholo_Yang--Mills_Proposition_}
follows from the

\hfill

\lemma \label{_Lambda_of_inva_forms_zero_Lemma_}
Let $\Theta\in \Lambda^2(M)$ be a $SU(2)$-invariant 
differential 2-form on $M$. Then
$\Lambda_L(\Theta)=0$ for each induced complex structure
$L$.\footnote{By $\Lambda_L$ we understand the Hodge operator 
$\Lambda$ associated with the K\"ahler complex structure $L$.}

{\bf Proof:} This is Lemma 2.1 of \cite{_Verbitsky:Hyperholo_bundles_}.
\endproof
 
\hfill

Let $M$ be a compact hyperk\"ahler manifold, $I$ an induced 
complex structure. 
For any stable holomorphic bundle on $(M, I)$ there exists a unique
Hermitian Yang-Mills connection which, for some bundles,
turns out to be hyperholomorphic. It is possible to tell when
this happens (though in the present paper we never use
this knowledge).

\hfill

\theorem
Let $B$ be a stable holomorphic bundle over
$(M,I)$, where $M$ is a hyperk\"ahler manifold and $I$
is an induced complex structure over $M$. Then 
$B$ admits a compatible hyperholomorphic connection if and only
if the first two Chern classes $c_1(B)$ and $c_2(B)$ are 
$SU(2)$-invariant.\footnote{We use \ref{_SU(2)_commu_Laplace_Lemma_}
to speak of action of $SU(2)$ in cohomology of $M$.}

{\bf Proof:} This is Theorem 2.5 of
 \cite{_Verbitsky:Hyperholo_bundles_}. \endproof


\subsection{Generic holomorphically symplectic manifolds}


In this section, we define generic holomorphically symplectic manifolds.
Such manifolds, as seen later (\ref{_hyperkae_embeddings_Corollary_}),
admit a hyperk\"ahler structure $\c H$
such that every closed complex analytic subvariety is trianalytic
with respect to $\c H$
(for the definition of trianalytic subvarieities,
see \ref{_trianalytic_Definition_}). 
We follow \cite{Verbitsky:Symplectic_I_} 
(see also \cite{Verbitsky:Symplectic_II_}).

\hfill

Let $M$ be a compact holomorphically symplectic K\"ahler manifold. 
By \ref{_symplectic_=>_hyperkahler_Proposition_}, 
$M$ has a unique hyperk\"ahler structure with
the same K\"ahler class and holomorphic symplectic form. 
Therefore one can without ambiguity speak
about the action of $SU(2)$ on $H^*(M,\R)$ (see 
\ref{_SU(2)_commu_Laplace_Lemma_}).
Of course, this action essentially depends on the choice
of K\"ahler class.

\hfill

\definition \label{_generic_manifolds_Definition_} 
Let $\omega\in H^{1,1}(M)$ be the K\"ahler class  
of a holomorphically symplectic
manifold $M$. We say that $\omega$ {\bf induces the $SU(2)$-action
of general type} when all elements of the group
\[ \bigoplus\limits_p H^{p,p}(M)\cap H^{2p}(M,\Z)\subset H^*(M)\] 
are $SU(2)$-invariant. 
A holomorphically symplectic manifold $M$ is 
called {\bf of general type} if there
exists a K\"ahler class on $M$ which induces 
an $SU(2)$-action of general type.

\hfill

As \ref{_gene_type_co_div_by2_Remark_} below
implies, holomorphically symplectic manifolds
of general type have no Weil divisors.  
In particular, such manifolds are never algebraic.

\hfill

\proposition \label{_generic_are_dense_Proposition_} 
Let $M$ be a hyperk\"ahler manifold and $S$
be the set of induced complex structures over $M$. Denote by 
$S_0\subset S$ the set of $L\in S$ such that the natural 
K\"ahler metric on $(M,L)$ induces the $SU(2)$ action of 
general type. Then $S_0$ is dense in $S$.

{\bf Proof:} This is Proposition 2.2 from
\cite{Verbitsky:Symplectic_II_}
\endproof

\hfill

One can easily deduce from the results in 
\cite{_Todorov:Moduli_I_II_} and
\ref{_generic_are_dense_Proposition_} that the set of points 
associated with holomorphically symplectic 
manifolds of general type is dense in the classifying space
of holomorphically symplectic manifolds of K\"ahler type.


\subsection{Trianalytic subvarieties in compact hyperk\"ahler
manifolds.}


In this subsection, we give a definition and a few basic properties
of trianalytic subvarieties of hyperk\"ahler manifolds. 
We follow \cite{Verbitsky:Symplectic_II_}.

\hfill

Let $M$ be a compact hyperk\"ahler manifold, $\dim_\R M =2m$.

\hfill

\definition\label{_trianalytic_Definition_} 
Let $N\subset M$ be a closed subset of $M$. Then $N$ is
called {\bf trianalytic} if $N$ is a complex analytic subset 
of $(M,L)$ for any induced complex structure $L$.

\hfill

Throughout this paper, we implicitly assume that our
trianalytic subvarieties are connected. Most results are trivially
generalized to the general case.

\hfill

Let $I$ be an induced complex structure on $M$,
and $N\subset(M,I)$ be
a closed analytic subvariety of $(M,I)$, $dim_\C N= n$.
Denote by $[N]\in H_{2n}(M)$ the homology class 
represented by $N$. Let $\inangles N\in H^{2m-2n}(M)$ denote 
the Poincare dual cohomology class. Recall that
the hyperk\"ahler structure induces the action of 
the group $SU(2)$ on the space $H^{2m-2n}(M)$.

\hfill

\theorem\label{_G_M_invariant_implies_trianalytic_Theorem_} 
Assume that $\inangles N\in  H^{2m-2n}(M)$ is invariant with respect
to the action of $SU(2)$ on $H^{2m-2n}(M)$. Then $N$ is trianalytic.

{\bf Proof:} This is Theorem 4.1 of 
\cite{Verbitsky:Symplectic_II_} (see Subsection
\ref{_SU(2)-inv=>triana_Subsection_} for a sketch of a proof). 
\endproof

\remark \label{_triana_dim_div_4_Remark_}
Trianalytic subvarieties have an action of quaternion algebra in
the tangent bundle. In particular,
the real dimension of such subvarieties is divisible by 4.

\hfill

\ref{_G_M_invariant_implies_trianalytic_Theorem_} has the following
immediate corollary:

\corollary \label{_hyperkae_embeddings_Corollary_} 
Let $M$ be a compact holomorphically symplectic 
manifold of general type, $S\subset M$ be its complex analytic
subvariety, and $\omega$ be a K\"ahler class
which induces an $SU(2)$-action of general type.
Consider the hyperk\"ahler structure associated with 
$\omega$ by \ref{_symplectic_=>_hyperkahler_Proposition_}.
Then $S$ is trianalytic with respect to $\c H$.

\endproof

\remark \label{_gene_type_co_div_by2_Remark_}
{}From \ref{_hyperkae_embeddings_Corollary_} and
\ref{_triana_dim_div_4_Remark_}, it follows that
a holomorphically symplectic manifold of general type
has no closed complex analytic subvarieties of odd dimension;
in particular, such a manifold has no divisors.


\subsection{Wirtinger's inequality and its use in K\"ahler geometry}


Let $M$ be a compact K\"ahler manifold, $X\subset M$ a 
closed real analytic subvariety.
In this section, we give criteria for $X$ to be 
complex analytic, in terms of certain integrals.
We follow \cite{Verbitsky:Symplectic_II_} and 
\cite{_Stolzenberg_}.

\hfill

\definition \label{_volume_Riema_Definition_}
Let $H$ be an $\R$-linear space equipped with a
positively defined scalar product, $\dim H =h$. The exterrior
form $\Vol\in \Lambda^h(H)$ is called {\bf a volume form}
if the the standard hypercube with the side 1
has the volume 1 in the measure defined by $\Vol$.

\hfill

Clearly, the volume form is defined up to a sign.
This sign is determined by the choice of orientation on $H$.
In the same manner 
we define the top degree differential form $\Vol$ called
{\bf a volume form} on any oriented 
Riemannian manifold. 

\hfill

Let $V$ be a Hermitian linear space, $W\subset V$ be a 
$\R$-linear subspace, $\dim_\R W = 2n$. Consider 
space $\Lambda^{2n}(W)$ of volume forms on $M$. Let
$\omega$ be the imaginary part of the Hermitian form 
on $V$. Consider the vectors $\omega^n$, 
$\Vol \in \Lambda^{2n}(W)$. Since $\Vol$ is non-zero,
and $\Lambda^{2n}(W)$ is one-dimensional, we can speak
of a fraction $\frac{\omega^n}{\Vol}$, which is a real number,
defined up to a sign (the form $\Vol$ is defined up to a sign).
Denote by $\Xi_W$ the number 
$\left|\frac{\omega^n}{\Vol}\right|$.

\hfill

\lemma\label{_Wirtinger_Lemma_} 
(Wirtinger's inequality) 
In these assumptions, $\Xi_W\leq 2^n.$
Moreover, if $\Xi_W=2^n$, then $W$ is a 
complex subspace of $V$.

{\bf Proof:} \cite{_Stolzenberg_} page 7. \endproof

\hfill

Let $M$ be a K\"ahler manifold, $N\subset (M,I)$ be a closed
real analytic subvariety of even dimension, 
$N_{ns}\subset N$ the set of non-singular 
points of $N$. For each
$x\in N_{ns}$, consider $T_xN$ as a subspace of $T_xM$.
Denote the corresponding number
$\Xi_{T_xN}$ by $\Xi_L (x)$. The following statement is a direct
consequence of \ref{_Wirtinger_Lemma_}:

\hfill

\proposition\label{_N_is_analytic_if_eta_is_constant_Proposition_} 
Let $J$ be an induced complex structure.
The subset $N\subset M$ is complex analytic 
if and only if 

\[ \forall x\in N_{ns} \;\;\;\; \Xi_L(x)=2^n. \]

\endproof

\hfill

Let $\Vol N_{ns}$ be the volume form of $N_{ns}$, taken with respect to the
Riemannian form (see \ref{_volume_Riema_Definition_}).

\hfill

\theorem \label{_Wirti_for_Kahle_Theorem_}
Let $M$ be a K\"ahler manifold, $N\subset M$, $\dim N =2n$ a 
closed real analytic subvariety, $N_{ns}\subset M$ the set
of its non-singular points. Assume that the improper integrals
$\int_{N_{ns}} \Vol N_{ns}$, $\int_{N_{ns}} \omega^n$ exist.
Then
\[ 
   2^n\int_{N_{ns}} \Vol N_{ns} \geq \int_{N_{ns}} \omega^n 
\]
and the equality is reached if and only if $N$ is complex analytic
in $M$.

{\bf Proof:} This is a direct consequence of 
\ref{_N_is_analytic_if_eta_is_constant_Proposition_} 
\endproof

\hfill

We use the term ``symplectic volume'' for the number
$\frac{1}{2^n}\int_{N_{ns}} \omega^n $ and ``Riemannian volume''
for $\int_{N_ns} \Vol N_{ns}$. Then, 
\ref{_Wirti_for_Kahle_Theorem_} might be rephrased in the form
``a real analytic cycle is complex analytic if and only if
its symplectic volume is equal to its Riemannian volume''.


\subsection{$SU(2)$-invariant cycles in cohomology 
and trianalytic subvarieties}
\label{_SU(2)-inv=>triana_Subsection_}


This subsection consists of a sketch of a proof of 
\ref{_G_M_invariant_implies_trianalytic_Theorem_}.
We follow \cite{Verbitsky:Symplectic_II_}.

\hfill

For a K\"ahler manifold $M$, $m=dim_\C M$ 
and a form $\alpha\in H^{2p}(M,\C)$, define 
\[
  \deg(\alpha):=\int_M \omega^{m-p}\wedge\alpha,
\]
where $\omega$ is the K\"ahler form.

Let $M$ be a compact hyperk\"ahler manifold, $\alpha$ a
differential form.
We denote by $\deg_L\alpha$ the degree associated with
an induced complex structure $L$.

\hfill

\proposition \label{_G_M_invariant_cycles_over_Proposition_} 
Let $M$ be a compact hyperk\"ahler manifold and
$\alpha$ be an $SU(2)$-invariant form of non-zero degree.
Then the dimension of $\alpha$ is divisible by 4. Moreover, 
\[ 
   \deg_{I}\alpha = \deg_{I'}\alpha,
\]
for every pair of induced complex structures $I$, $I'$.

{\bf Proof:} This is Proposition 4.5 of 
\cite{Verbitsky:Symplectic_II_}.
\endproof

\hfill

\ref{_G_M_invariant_implies_trianalytic_Theorem_} follows
immediately from \ref{_G_M_invariant_cycles_over_Proposition_} 
and \ref{_Wirti_for_Kahle_Theorem_}. By
\ref{_Wirti_for_Kahle_Theorem_}, a real analytic subvariety
$N \subset M$ is trianalytic if and only if its symplectic volume,
taken with respect to any of induced complex structures, is 
defined and equal to its Riemannian volume.
In notation of \ref{_G_M_invariant_implies_trianalytic_Theorem_},
the symplectic volume of $N$ taken with respect to $L$
is equal to $\frac {1}{2^n} \deg_L\inangles N$.
{}From \ref{_G_M_invariant_cycles_over_Proposition_},
we obtain that symplectic volume of $N$ is the same for all
induced complex structures. Since $N$ is complex analytic
with respect to $I$, the symplectic volume of $N$ is equal
to the Riemannian volume of $N$, again by 
\ref{_Wirti_for_Kahle_Theorem_}. This proves
\ref{_G_M_invariant_implies_trianalytic_Theorem_}.


\section[Deformations of non-singular trianalytic subvarieties.]
{Deformations of non-singular trianalytic \\subvarieties.}
\label{_nonsingu_preli_Section_}


Let $M$ be a holomorphically symplectic K\"ahler manifold,
$X\subset M$ a closed complex analytic submanifold 
of trianalytic type. We study
the deformations of $X$ in $M$, in order to prove 
\ref{_iso_intro:Theorem_}. It turns out that
\ref{_iso_intro:Theorem_} is almost trivial in the case
$X$ non-singular. In this section, we give the proof
of \ref{_iso_intro:Theorem_} (i) for non-singular $X$; the general
case is proven independently. We hope that a simple argument
will be insightful, even if we need to produce a separate proof
for the general case. This section is perfectly safe to skip.

\hfill

\remark For a smooth trianalytic submanifold $X\subset M$,
$X$ is obviously hyperk\"ahler, in a natural way. The quaternion
action comes from quaternion action in $TM$, and the metric is
induced from $M$ too.

\hfill

We start from the following simple, but important, lemma.

\hfill

\lemma \label{_TM-restrict_X_hyperholo:Lemma_} 
Let $M$ be a complex manifold equipped with a hyperk\"ahler metric,
$X\subset M$ a closed complex analytic submanifld of trianalytic
type. Consider the
restriction $TM\restrict{X}$ of the tangent bundle to $M$ on $X$.
We equip $TM_\restrict{X}$ with a connection $\nabla$ coming
from Levi-Civita connection in $TM$. Then $(TM\restrict{X}, \nabla)$
is {\bf hyperholomorphic}.\footnote{See 
Subsection \ref{_hyperholo_Subsection_} 
for the definition of {\bf hyperholomorphic}.
By \ref{_hyperholo_Yang--Mills_Proposition_},
a hyperholomorphic connection  is Yang-Mills.}

\hfill

{\bf Proof:} Consider the natural action of $SU(2)$ on the 
space 
\[ 
   \Lambda^2\left(X, \End\left(TM\restrict{X}\right)\right).
\] 
We need to show that
the curvature $\Theta$ of $\nabla$ is $SU(2)$-invariant. Let
$\nabla_{LC}$ be the Levi--Civita connection on $TM$, and 
$(\nabla_{LC})^2\in \Lambda^2(M, \End(TM))$ be its curvature. Clearly,
$\Theta$ is a pull-back of $(\nabla_{LC})^2$ to $X$. Therefore,
it suffices to show that $(\nabla_{LC})^2$ is $SU(2)$-invariant.
This is \ref{_tangent_hyperholo_Example_}. 
\endproof

\hfill

As a corollary, we obtain the following proposition.

\hfill

\proposition \label{_TM-restrict-X_decompo_for_compa_Proposition_}
In assumptions of \ref{_TM-restrict_X_hyperholo:Lemma_},
let $M$, $X$ be compact. Then the following statements are true.

\begin{description}
\item[(i)] The bundle $TM\restrict{X}$ is naturally 
isomorphic to the direct sum $TM\restrict X \cong TX \oplus NX$,
where $N$ is the normal bundle to $X$ inside $M$. This
isomorphism is compatible with the natural connections and Hermitian
metrics on $TM\restrict X$, $TX$, $NX$. 

\item[(ii)] $NX$ is hyperholomorphic.

\item[(iii)] For each section $\gamma$ of $NX$, $\gamma$ is nowhere
degenerate, and there is a natural decomposition
\begin{equation} \label{_NX=O-gamma+rest:Equation_}
NX = \calo_\gamma \oplus N_\gamma X,
\end{equation}
where $\calo_\gamma$ is a trivial sub-bundle of $NX$ generated by
$\gamma$, and $N_\gamma X$ its orthogonal complement. The 
decomposition \eqref{_NX=O-gamma+rest:Equation_} is compatible 
with connection. 

\item[(iv)] In assumptions of (iii),
consider the connection $\nabla_\gamma$ induced from
$NX$ to $\calo_\gamma X$. Then is $\nabla_\gamma$ flat.
\end{description}

\hfill

\remark In fact, 
\ref{_TM-restrict-X_decompo_for_compa_Proposition_} (i)-(ii) holds true 
in assumptions of \ref{_TM-restrict_X_hyperholo:Lemma_} for
general, not necessarily compact, $X$ and $M$ 
(see \ref{_NX_splits_for_hype_Proposition_}, 
\ref{_NX_hyperholo_Proposition_}).
However, the proof is easier in compact case.

\hfill

{\bf Proof of \ref{_TM-restrict-X_decompo_for_compa_Proposition_}.}
Consider the embedding 
\[ TX \hookrightarrow TM\restrict X \]
of bundles with connection. By \ref{_tangent_hyperholo_Example_},
$TX$ is hyperholomorphic, and hence Yang--Mills.
By \ref{_TM-restrict_X_hyperholo:Lemma_}, $TM\restrict X$ is
hyperholomorphic as well. We obtain that
$TX$ is a destabilizing subsheaf in $TX$.
By \ref{_UY_Theorem_}, a Yang--Mills bundle 
is a direct sum of stable bundles.\footnote{Bundles which are direct
sum of stable are called {\bf polystable}.}
Thus, $TM$ is a direct sum of $TX$ and its orthogonal
complement $NX$ (see \ref{_YM_exact_split_Proposition_} 
for details). This proves
\ref{_TM-restrict-X_decompo_for_compa_Proposition_} (i).
Since the curvature
of $TM\restrict X$ is decomposed onto a 
direct sum of the curvature of $NX$ and the curvature of
$TX$, the curvature of $NX$ is $SU(2)$-invariant.
This proves \ref{_TM-restrict-X_decompo_for_compa_Proposition_} (ii).

\hfill

The condition (iii) follows from (ii) and \ref{_UY_Theorem_}, 
since every holomorphic section $\gamma$ of a 
Yang--Mills bundle $B$ of zero degree (such as hyperholomorphic bundles)
spans a destabilizing subsheaf $\calo_\gamma$. 
Since $B$ is Yang--Mills, it is polystable; thus, 
$\calo_\gamma$ is its direct summand.
This proves \ref{_TM-restrict-X_decompo_for_compa_Proposition_} 
(iii). To prove (iv), we notice that $\calo_\gamma$
is a trivial holomorphic bundle with Yang--Mills connection.
By \ref{_UY_Theorem_}, Yang--Mills connection is unique,
and therefore, $\calo_\gamma$ is flat. 
\ref{_TM-restrict-X_decompo_for_compa_Proposition_} 
is proven. 
\endproof

\hfill

We recall the following general results of the
theory of deformations of complex subvarieties.
Let $X\subset M$ be a closed complex analytic subvariety
of a compact complex manifold. The {\bf Douady space}
of deformations of $X$ inside of $M$ is defined
(\cite{_Douady_}). We denote the Douady space by $D_M X$,
or sometimes by $D(X)$. By definition, $D_M X$ is a
finite-dimensional complex analytic variety.
The points of $D_MX$ are identified with the 
subvarieties $X'\subset M$, where $X'$ is a deformation
of $X$. Let $\gamma:\; \R \arrow D_m X$ be a real analytic map.
Assume that for $t=t_0$, $\gamma(t)$ is a smooth subvariety of
$M$. In deformation theory, the differential
$\frac{d\gamma}{dt}(t_0)$ is interpreted as a 
holomorphic section of the normal bundle $N \gamma (t_0)$.
Thus, the Zariski tangent space $T_{X} D_M (X)$
is naturally embedded to the space of holomorphic
sections of $N X$. The following proposition is an
easy application of the Kodaira--Spencer theory.

\hfill

\proposition \label{_compa_subva_iso_holo:Proposition_}
Let $X \subset M$ be a closed complex analytic submanifold of
a compact complex manifold $M$, and $D_M X$ the corresponding
Douady space. Let $\gamma:\; \R \arrow D_M X$ be a real analytic
map. Assume that for all $t_0\in \R$, the subvariety
$\gamma(t_0)$ is smooth. Assume also that the section 
\[ 
  \nu_{t_0} = 
  \frac{d \gamma}{dt}(t_0) \in \Gamma_{\gamma(t_0)} (N \gamma(t_0))
\]
is nowhere degenerate. Finally, assume that $\nu_{t_0}$
splits from $N \gamma(t_0)$, i. e, there exists a decomposition
of a holomorphic vector bundle
\[
   N \gamma (t_0) = \calo_{\nu_{t_0}} \oplus N'
\]
where $\calo_{\nu_{t_0}}$ is the trivial subbundle of 
$N \gamma(t_0)$ generated by $\nu_{t_0}$.
Then the subvarieties $\gamma(t)\subset M$ are naturally
isomorphic for all $t$.

{\bf Proof:} Well known; see, for instance,
\cite{_Kodaira_Spencer_} \endproof

\hfill

\proposition \label{_compa_subva_iso_metric:Proposition_}
In assumptions of \ref{_compa_subva_iso_holo:Proposition_},
let $M$ be K\"ahler. Consider the induced 
Hermitian structure on $N\gamma(t)$, for all $t\in \R$.
Assume that for all $t\in \R$, the section
\[ 
  \nu_{t_0} = 
  \frac{d \gamma}{dt}(t_0)
\]
has constant length. Assume also that the 
orthogonal decomposition \[ TM\restrict X  = NX \oplus NX^\bot \]
is compatible with the Levi-Civita connection in 
$TM\restrict X$. Consider the isomorphisms 
$\psi_{t_1,t_2}:\; \gamma(t_1)\arrow \gamma(t_2)$ constructed in
\ref{_compa_subva_iso_holo:Proposition_}. Then 
the maps $\psi_{t_1,t_2}$ are compatible with the K\"ahler
metric induced from $M$.

{\bf Proof:} The proof follows from Kodaira--Spencer construction;
for a complete argument, see 
\ref{_conne_in_fami_of_comple_geo_Proposition_} (ii). 
\endproof

\hfill

The following theorem is the main result of this section.

\hfill

\theorem \label{_iso_for_smooth_subva:Theorem_}
Let $M$ be a compact holomorphically symplectic 
K\"ahler manifold,
and $X\subset M$ a complex submanifold
of trianalytic type. Let $X'$ be a deformation of $X$ 
in $M$. Then there exists a complex analytic
isomorphism $\psi:\; X \arrow X'$. Moreover, if
the K\"ahler metric on $M$ is hyperk\"ahler,
then $\psi:\; X \arrow X'$
is compatible with the K\"ahler structure induced from $M$.

\hfill

{\bf Proof:} The deformations of $X$ are infinitesimally classified
by the sections of $NX$. Let $\gamma$ be such a section.
Applying \ref{_TM-restrict-X_decompo_for_compa_Proposition_},
we obtain that assumptions of
\ref{_compa_subva_iso_holo:Proposition_} and
\ref{_compa_subva_iso_metric:Proposition_} hold for the deformations
of $X$ inside $M$. Now, \ref{_iso_for_smooth_subva:Theorem_}
is directly implied by \ref{_compa_subva_iso_holo:Proposition_} and
\ref{_compa_subva_iso_metric:Proposition_}. \endproof

\hfill

\remark 
\ref{_iso_for_smooth_subva:Theorem_} 
gives a proof of \ref{_iso_intro:Theorem_} 
(i), for the subvarieties $X\subset M$ which are
smooth. We prove \ref{_iso_for_smooth_subva:Theorem_} for
general subvarieties in \ref{_triana_subse_comple_ana_Theorem_}.


\section[Completely geodesic embeddings of Riemannian manifolds.]
{Completely geodesic embeddings \\of Riemannian manifolds.}
\label{_comple_geode_defo_Section_}


Before we proceed with the proof of \ref{_iso_intro:Theorem_},
we have to prove a serie of preliminary 
results from deformation theory. The arguments of
deformation theory are greatly simplified in the
hyperk\"ahler case, because hyperk\"ahler embeddings are
completely geodesic\footnote{For a definition of completely geodesical
embeddings, see \ref{_comple_geode:Definition_}.}
 (\ref{_hype_embe_comple_geode:Corollary_}).
We prove a number of simple statements from the deformation
theory of completely geodesical embeddings of K\"ahler manifolds.

\subsection{Completely geodesic submanifolds}

\nopagebreak
\hspace{5mm}
\proposition \label{_comple_geodesi_basi_Proposition_}
Let $X \stackrel \phi\hookrightarrow M$ be an embedding of Riemannian 
manifolds (not necessarily compact) compatible with the Riemannian
structure.\footnote{Such embeddings are called 
{\bf Riemannian embeddings}.}
 Then the following conditions are equivalent.

\begin{description}
\item[(i)] For every point $x\in X$, there exist a neighbourhod $U \ni x$
of $x$ in $X$ such that for all 
$x' \in U$ there is a geodesic in $M$
going from $\phi(x)$ to $\phi(x')$ which lies 
in $\phi(X)\subset M$.

\item[(ii)] Consider the Levi-Civita connection $\nabla$ on $TM$,
and restriction of $\nabla$ to $TM \restrict{X}$. Consider the
orthogonal decomposition 
\begin{equation} \label{TM_decompo_Equation_} 
   TM\restrict{X} = TX \oplus TX^\bot. 
\end{equation}
Then, this decomposition is preserved by the connection $\nabla$.
\end{description}

{\bf Proof:} Well known; see, for instance, 
\cite{_Besse:Einst_Manifo_}.
\blacksquare

\hfill

\definition\label{_comple_geode:Definition_}
Let $X \stackrel i \hookrightarrow M$ be a Riemannian embedding
satisfying either of the conditions of 
\ref{_comple_geodesi_basi_Proposition_}. Then $i$ is called 
{\bf a completely geodesic embedding}, and the image $i(X)\subset M$
is called {\bf a completely geodesic submanifold}.

\hfill

\lemma \label{_comple_geo_compa_holo_Lemma_}
Let $M$ be a K\"ahler manifold, $X\subset M$ a complex submanifold.
If $X$ is completely geodesic, then the decomposition
\eqref{TM_decompo_Equation_}  
is compatible with the holomorphic structure
on $TX$, $TM\restrict X$.

{\bf Proof:} Clear. \endproof


\subsection{Deformations of submanifolds}
\label{_defo_subva_conven_Subsection_}


Let $\c X \stackrel \pi \arrow S$ be a family of complex
manifolds equipped with a map $\c X \stackrel \phi \arrow M$.
Denote the pre-image $\pi^{-1}(s)\subset \c X$ by $X_s$.
Assume that for all $s\in S$
the restriction $\phi_s:\; X_s\arrow M$ of 
$\c X \stackrel \phi \arrow M$
to $X_s \subset \c X$ is a smooth embedding.
Assume that for $s_0\in S$, the image $\phi_{s_0}(X_{s_0})\subset M$
coinsides with $X$.

\hfill

\definition 
The collection of data 
\[ \left( \c X \stackrel \pi \arrow S, 
   \phi, \phi_{s_0}:\; X_{s_0} \oldtilde \arrow X\right )
\]
is called {\bf a family of submanifolds of $M$},
and {\bf a family of deformations of $X$}.
The same definition can be formulated for $M$, $X$, $S$, $\c X$
real analytic; in this case, we speak of {\bf family of real
analytic submanifolds of $M$}. Also, the submanifolds might be 
replaced by subvarieties in order to obtain the definition
of a family of subvarieties.


\subsection{Section of a normal bundle arising from a deformation}
\label{_norma_vecto_Subsection_}


Let $M$ be a complex or real analytic manifold and
\[ \left( \c X \stackrel \pi \arrow S, 
   \phi:\; \c X \arrow M \right)
\]
be a system of submanifolds.
For each tangent vector $t \in T_{s_0} S$, deformation theory 
gives a canonical section $\eta$ of the normal bundle $N X_{s_0}$.
This section is holomorphic if we work in complex situation.
We recall how this section is obtained.

For a sufficiently small neighbourhood of $x\in \c X$,
we can always find coordinates in $X_s$ which analytically
depend on $s$. For each point $x\in X_{s_0}$, 
denote by $x(s)$ the point of $X_s$ with the same
coordinates as $x$. Consider the vector
$\tilde \eta_x :=\frac{dx(s)}{ds}(t)\in T_xM$, which is a derivative
of $x(s)$ along $t$. The vector $\tilde \eta_x$ obviously depends
on the choice of 
coordinates in $X_s$. Let $\eta_x\in NX_{s_0}\restrict x$ be 
the image of $\tilde\eta_x$ under the natural map 
$T_x M \arrow N X_{s_0}\restrict x$. Clearly, $\eta_x$ is
independent from the choice of coordinates. Gluing
$\eta_x$ together, we obtain the canonical 
section $\eta\in NX_{s_0}$.

We need the following technical lemma in Section 
\ref{_triholo_Section_}.

\hfill

\lemma \label{_norma_sec_holom_Lemma_}
Let $M$ be a complex analytic manifold and
\[ \left( \c X \stackrel \pi \arrow S, 
   \phi:\; \c X \arrow M \right)
\]
be a real analytic system of submanifolds. Assume that for 
all $s\in S$, the submanifold $\pi^{-1}(s) := X_s\subset M$ 
is complex analytic.
Then, for all $t\in T_{s_0} S$, the corresponding section
$\eta \in N X_{s_0}$ is holomorphic.

\hfill

{\bf Proof:} Shrinking $\c X$ if necessary, we 
can find a system of holomorphic
coordinates in $X_s$, $s\in S$ 
which depends smoothly upon $s$. Consider the section
$\tilde \eta\in TM\restrict{X_{s_0}}$ obtained by
deriving $x(s)$ along $t$ as above. Since the map
$x\arrow x(s)$ is by construction holomorphic,
the section $\tilde\eta$ is also holomorphic. Then,
$\eta$ is holomorphic by construction. \endproof


\subsection{Connections in families of manifolds}
\label{_conne_Subsection_}


\hfill

\definition \label{_conne_Definition_}
Let $\pi:\; \c X \arrow S$ be a family
of real analytic manifolds parametrized by $S$. 
Consider the bundles $T\c X$, $N_\pi X$,
where $N_\pi X$ is a fiberwise normal bundle 
to the fibers of $\pi$. There is a natural 
projection $p:\; T\c X \arrow N_\pi X$.
The {\bf connection} in $\c S$ is 
a section of $p$, i. e. such an embedding
$\nabla:\; N_\pi X \arrow T\c X$ that 
$p \circ s = \Id \restrict {N_\pi X}$. 

This definition
is naturally adopted to the case of complex analytic family
of manifolds. If the 
section $\nabla:\; N_\pi X \arrow T\c X$ is holomorphic,
the connection is called {\bf a holomorphic connection}. 

This definition makes sense only when the base variety $S$ is smooth.
However, the connection in a family gives a connection on a pullback
of this family, under all maps $s:\; S' \arrow S$. Connections on the 
pullback families are naturally compatible. We shall sometimes speak 
of {\bf connection} in a family where the base $S$ is not smooth.
This means that for all maps $s:\; S' \arrow S$, where $S'$ is smooth,
the pullback family is equipped with a connection, and these
connections are compatible. All the definitions and results
(which we state and prove in the case of a smooth base) are
naturally adopted to the case when the base $S$ is singular.

\hfill

\proposition \label{_conne_in_fam_of_comple_geode_Proposition_}
Let 
\[ \left( \c X \stackrel \pi \arrow S, 
   \phi:\; \c X \arrow M \right)
\]
be a family of submanifolds of a Riemannian manifold $M$.
Assume that for all $s\in S$, the submanifold 
$X_s:= \psi(\pi^{-1}(s)) \subset M$ is completely geodesic. Then
the family $\c X$ is equipped with a natural connection.
Moreover, if $M$ is K\"ahler and the family $\c X$ is complex
analytic, then the connection $\nabla$ is holomorphic.

\hfill

{\bf Proof:} Let $T_\pi X\subset T\c X$ be the bundle of vectors
tangent to the fibers of $\pi$.
To split the exact sequence
\[ 
    0 \arrow T_\pi X \stackrel i\arrow T\c X \arrow T_\pi X \arrow 0,
\]
we have to construct a surjection 
\begin{equation} \label{_surje_secti_for_conne_Equation_}
   T{\c X} \stackrel p \arrow T_\pi X
\end{equation}
satisfying $i\circ p = \Id_{T_\pi X}$.
Consider the orthogonal decomposition 
\begin{equation} \label{_ortho_compo_for_norma_in_fami_Equation_}
\phi^*TM=\phi^*N_M X_t \oplus TX_t,
\end{equation}
where $N_M X_t$ is the normal bundle to $\phi(X_t)$ in $M$.
This gives a natural epimorphism 
\[ \phi^* TM \stackrel e\arrow T_\pi X.\] Taking a composition
of $e$ with \[ d\phi:\; T\c X \arrow \phi^* TM,\]
we obtain a surjection $p$ of 
\ref{_surje_secti_for_conne_Equation_}. 

In K\"ahler case, the decomposition
\eqref{_ortho_compo_for_norma_in_fami_Equation_}
is holomorphic by \ref{_comple_geo_compa_holo_Lemma_}. Thus, 
we constructed a connection which is holomorphic.
\endproof

\hfill

For each real analytic path $\gamma:\; [0,1] \arrow S$,
$t, t'\in [0,1]$, we may {\bf integrate} the connection
along $\gamma$. This notion is intuitively
clear (and well known). We recall briefly
the definition of the integral maps associated to connection,
in order to fix the notation.

Restricting the family $\c X$ to the image of $\gamma$
in $S$ does not change the result of integration.
To simplify the exposition, we assume that $S$ coinsides with
the image $D \cong [0,1]$ of $\gamma$. In this case, the bundle $N_\pi X$ 
is naturally isomorphic to a pullback $\pi^* T D$. Denote by
$\nu$ a canonical unit section of $\pi^* TD$ corresponding to
the unit tangent field to $[0,1]$. Then, $\nabla \nu$ is a vector
field in $T\c X$. Integrating $\nu$ to a 
diffeomorphism, we obtain a map $\exp(t \nu)$ defined in an open
subset of $\c X$, for $t\in \R$. Clearly, $\exp(t\nu)$ maps the points
of $X_{t_1}$ to the points of $X_{t_1+t}$. The resulting
diffeomorphism we denote by 
\[ 
   \Psi^{t_1}_{t_1+t}:\; U_{t_1} \arrow U_{t+t_1}, 
\]
where $U_{t_1}$ is an intersection of $X_{t_1}$ with the domain 
of $\exp(t\nu)$. Clearly, the same definition applies to the
arbitrary families of manifolds with connection. When $\c X$
is a complex analytic family with a holomorphic connection,
the diffeomorphisms $\Psi^t_{t'}$ are complex analytic.


\subsection{Deformations of completely
geodesic submanifolds (the main statement)}


Let  
\[ \left( \c X \stackrel \pi \arrow S, 
   \phi, \phi_{s_0}:\; X_{s_0} \oldtilde \arrow X\right )
\]
be a real analytic deformation of $X\subset M$,
and $\gamma:\; [0,1] \arrow S$ be a real analytic map.
Slightly abusing the notation, we 
shall consider the fibers $X_s$ of $\pi$ as subvarieties in $M$.
Assume that $X_{\gamma(t)}$ is a completely geodesic complex
analytic submanifold of $M$, for all
$t$. Let $\eta_t \in \Gamma \left(N X_{\gamma(t)}\right)$ be a section of
the normal bundle to $X_{\gamma(t)}$ corresponding to the vector 
$\frac {d \gamma}{dt} \in T_{\gamma(t)} S$. 
Using the natural 
holomorphic embedding 
\begin{equation} \label{_N_embe_to_TM_Equation_} 
  N X_{\gamma(t)} \hookrightarrow TM\restrict{X_{\gamma(t)}}
\end{equation}
associated with a decomposition 
$TM\restrict{X_{\gamma(t)}} = N X_{\gamma(t)}\oplus TX_{\gamma(t)}$,
we may consider $\eta_t$ as a section of
$TM\restrict{X_{\gamma(t)}}$.

\hfill

\proposition \label{_conne_in_fami_of_comple_geo_Proposition_}
Let $t$, $t'\in [0,1]$. 
Let $\Psi^t_{t'}:\; U_t \arrow U_{t'}$ be the map of
Subsection \ref{_conne_Subsection_} 
integrating the connection $\nabla$ of 
\ref{_conne_in_fam_of_comple_geode_Proposition_},
where $U_t$, $U_{t'}$ are open subsets 
of $X_{\gamma(t)}$, $X_{\gamma(t')}$ defined in 
Subsection \ref{_conne_Subsection_}.
Assume that for all $t$ the section 
$\eta_t\in TM\restrict{X_{\gamma(t)}}$ is 
parallel with respect to the natural
connection on $TM\restrict{X_{\gamma(t)}}$
obtained from Levi--Civita on $TM$. Assume also that the Riemannian form
on $M$ is real analytic.\footnote{This assumption is extraneous;
we use it to simplify the exposition. For the case we are interested
in ($M$ a hyperk\"ahler manifold) the Riemannian form
is real analytic by \ref{_Riema_on_hype_real_ana_Corollary_}.}
Then $\Psi^t_{t'}$ satisfies the following conditions.

\begin{description}

\item[(i)] $\Psi^t_{t'}$ is compatible
with a Riemannian structure on $U_t$, $U_{t'}$. 

\item[(ii)] 
Assume also that $M$ admits a smooth Riemannian
compactification $\bar M$ such that the family $\c X$ 
admits a compactification in $\bar M$ (such is the case
when $M$ and $\c X$ are quasiprojective, and $\phi$ is
algebraic). Then the maps $\Psi^t_{t'}:\; U_t \arrow U_{t'}$ of (i) 
can be found in such a way that $U_t$ is the set 
of non-singular points of $X_{\gamma(t)}$. 

\item[(iii)] In assumptions of (ii), we can extend $\Psi^t_{t'}$
to an isomorphism of metric spaces 
\[ 
   \bar \Psi^t_{t'}:\; \bar X_{\gamma(t)} \arrow \bar X_{\gamma(t')}, 
\]
where $\bar X_{\gamma(t)}$, $\bar X_{\gamma'(t)}$ are
closures of $X_{\gamma(t)}$, $X_{\gamma'(t)}$ in $\bar M$, with induced metrics.

\end{description}

The next subsection is taken by the proof of 
\ref{_conne_in_fami_of_comple_geo_Proposition_}.


\subsection{Deformations of completely
geodesic submanifolds (the proofs)}


The statement (i) is sufficient to prove
in a small neighbourhood, by analytic continuation. Thus, we
may pick an open subset $U\in {\c X}$ in such a way that the restriction of
$\phi:\; \c X \arrow M$ to 
\[ 
   \c X':= \pi^{-1}(\gamma([0,1])\cap U
\]
is an embedding. Shrinking $\c X$ to $\c X'$ and pulling
from $M$ the Riemannian metric, we see that
\ref{_conne_in_fami_of_comple_geo_Proposition_} (i) is 
implied by the following lemma.

\hfill

\lemma \label{_integra_field_Killing_Lemma_}
Let $\c X \stackrel \pi\arrow [0,1]$ 
be a real analytic family of manifolds, equipped with a Riemannian
metrics. Let $T _\pi X\subset T\c X$ 
be the relative tangent bundle consisiting
of all vectors tangent to the fibers of $\pi$, and 
$N _\pi X= T \c X /T _\pi X$ be the normal bundle 
to the fibration $\pi$. Assume that the fibers 
$\pi^{-1} (t) \subset \c X$ are completely geodesic 
in $\c X$. Consider the connection
$\nabla:\; N _\pi X \arrow T\c X$ of 
\ref{_conne_in_fam_of_comple_geode_Proposition_}.
Let $\Psi^t_{t'}:\; U_t \arrow U_{t'}$ be the
maps obtained by integrating the connection as in 
Subsection \ref{_conne_Subsection_}
and $\eta_t \in N _\pi X$ be the normal fields arising
from the deformation theory. Assume that $\eta_t$'s
are parallel with respect to the natural connection on $N _\pi X$.
Then the maps $\Psi^t_{t'}$ are compatible with Riemannian metrics.

\hfill

{\bf Proof:} Let $\eta\in T \c X$ be the tangent vector field
obtained by gluing all $\eta_t$ together. Then $\Psi_t$ 
can be considered as integral map of this vector field. Thus,
to prove that $\Psi_t$ is compatible with the Riemannian
structure, we have to show that 
$\eta$ is {\bf a Killing vector field}.\footnote{A Killing field 
is a vector field which integrates to a diffeomorphism which is
compatible with a Riemannian metric.} Denote by
$\nabla_x:\; T\c X \arrow T\c X$ the 
action of covariant derivative along the vector field $x$. 
By \cite{_Besse:Einst_Manifo_},
Theorem 1.81, to prove that $\eta$ is Killing it 
suffices to prove that for all fields $a, b\in T\c X$,
we have
\begin{equation}\label{_Killing_Equation_} 
   (\nabla_a \eta, b) + (\nabla_b \eta,a) =0,
\end{equation}
where $(\cdot,\cdot)$ is the Riemannian form.
Take coordinates $(x_0, x_1,...x_n)$ 
on $\c X$ in such a way that $x_0$ comes from a
projection $\pi:\; \c X \arrow [0,1]$ and 
$x_1, ... x_n$ are coordinates along the fibers.
Let $\frac{d}{d x_i}$ be the corresponding vector fields.
To prove that $\eta$ is Killing it suffices to check
\eqref{_Killing_Equation_} 
for $a$, $b$ coordinate vector
fields. Since $\eta$ is parallel along $X_t$, we have
$\nabla_{\frac{d}{d x_i}} \eta=0$ for $i= 1,..., n$.
Thus, it suffices to prove \eqref{_Killing_Equation_} 
in case $a = \frac{d}{dx_0}$. For appropriate choice
of coordinates, $\frac{d}{d x_0}=\eta$; thus, 
\eqref{_Killing_Equation_} is implied by
\begin{equation}\label{_nabla_eta_coo_Equation_}
   \left(\nabla_\eta \eta,\frac{d}{d x_i}\right) =0, 
   \;\; i= 0,\dots n.
\end{equation}
Since Levi-Civita connection is compatible with the metrics,
we have 
\begin{equation}\label{_Levi_Civi_compa_Equation_}
   D_\eta \left(\eta, \frac{d}{d x_i}\right) =
   \left(\nabla_\eta \eta,\frac{d}{d x_i}\right) + 
   \left(\eta,\nabla_\eta\frac{d}{d x_i}\right),
\end{equation}
where $D_\eta$ is the usual (directional) derivative along $\eta$.
For $i=0$, this gives 
\[ D_\eta (\eta, \eta) =
   (\nabla_\eta \eta,\eta) + 
   (\eta,\nabla_\eta\eta)=0
\]
($\eta$ is parallel, and hence has constant length).
This proves \eqref{_nabla_eta_coo_Equation_} for $i=0$.
Since $\frac{d}{d x_i}$ are coordinate vector fields
they commute, and therefore, $\nabla_\eta\frac{d}{d x_i}=
-\nabla_{\frac{d}{d x_i}}\eta=0$ for $i>0$. Thus,
the equation \eqref{_Levi_Civi_compa_Equation_} for $i>0$ is 
reduced to 
\[
   D_\eta \left(\eta, \frac{d}{d x_i}\right) =
   \left(\nabla_\eta \eta,\frac{d}{d x_i}\right). 
\]
Since the vector fields
$\frac{d}{d x_i}$, $i>0$ are tangent to the fibration $\pi$, the function
$(\eta, \frac{d}{d x_i})$ is identically zero. This proves 
\eqref{_nabla_eta_coo_Equation_} for $i>0$. 
\ref{_integra_field_Killing_Lemma_} and consequently
\ref{_conne_in_fami_of_comple_geo_Proposition_}
(i) is proven. \endproof

\hfill

To prove
\ref{_conne_in_fami_of_comple_geo_Proposition_}
(ii) we have to show that the connection of 
\ref{_conne_in_fam_of_comple_geode_Proposition_} 
can be integrated for all smooth points of $\c X$.
We give a sketch of a simple geometric argument.
By (i), the maps $\Psi_t$ are isometries. Thus, the 
distance from the given point 
to the singular set of $X_t$ is invariant under 
the maps $\Psi_t$. To integrate the connection,
we write a tangent vector field which we
subsequently integrate. Since the normal field $\eta_t$
is parallel, this tangent field is uniformly bounded.
On a certain distance from the singular set $Sing(X_t)$, depending
on this uniform bound, the connection might be always 
integrated. In more precise terms,
for all $\epsilon>0$, all $x\in X_t$, with the distance
between $x$ and $Sing(X_t)$ no less than $\epsilon$, there exist
$\delta>0$ and a map 
$\Psi^t_{t+\delta}:\;X_t \arrow X_{t+\delta}$ defined
in a neighbourhood of $x$ which integrates the connection.
{}From this statement and \ref{_conne_in_fami_of_comple_geo_Proposition_} 
(i),  \ref{_conne_in_fami_of_comple_geo_Proposition_} (ii)
follows directly.

Finally, to prove \ref{_conne_in_fami_of_comple_geo_Proposition_} (iii)
we notice that $X_t$ is completely geodesic in $M$ for all $t$. Thus,
the completion of $X_t$ as a metric space coinsides with the
closure of $X_t$ in $\bar M$. Every isometry of metric spaces 
extends to a completion, and thus, $\Psi^t_{t'}$ extends to a closure of
$X_t$, $X_{t'}$ in $\bar M$.
\endproof


\section[Hyperholomorphic bundles and completely geodesical 
embeddings.]{Hyperholomorphic bundles \\and completely geodesical
embeddings.} 
\label{_comple_geode_hyperho_Section_}


\hfill

In this section, we prove that hyperk\"ahler embeddings are
completely geodesic. 

\subsection{Hyperholomorphic structure on the normal bundle}

Let $M$ be a hyperk\"ahler manifold, not necessary compact,
and $B$ a vector bundle.
Recall that in Section \ref{_basics_Section_}, we defined 
{\bf hyperholomorphic connections} in $B$ 
(\ref{_hyperho_conne_Definition_}). 
These are connections $\nabla$ such that the curvature 
$\Theta:=\nabla^2\in \Lambda^2(M,End(B))$ is an 
$SU(2)$-invariant 2-form with respect to the
natural $SU(2)$-action in $\Lambda^2(M)$. 
The hyperholomorphic connections are always Yang-Mills 
(\ref{_hyperholo_Yang--Mills_Proposition_}). The Levi--Civita connection
in the tangent bundle $TM$ is a prime example of a hyperholomorphic
connection (\ref{_tangent_hyperholo_Example_}). 

\hfill

\proposition \label{_NX_hyperholo_Proposition_}
Let $M$ be a hyperk\"ahler manifold, not necessary compact, and 
$X\subset M$ a trianalytic submanifold. Consider the normal
bundle $NX$, equipped with a connection $\nabla$ 
induced from the Levi--Civita
connection on $M$. Then $NX$ is hyperholomorphic.

\hfill

{\bf Proof:} Let $L$ be an induced complex structure $M$.
Consider the manifold $(X, L)$ as a complex submanifold of $(M, L)$. The
normal bundle $NX = N(X,L)$ has a natural holomorphic structure,
which is compatible with the connection $\nabla$. Therefore,
the curvature $\Theta\in \Lambda^2(X, End(NX))$ is of type
$(1,1)$ with respect to the Hodge decomposition
defined by $L$. Thus, $\Theta$ is of type $(1,1)$ with
respect to any of the induced complex structures.
By definition, this means that $\nabla$ is hyperholomorphic.
\endproof

\subsection{Hyperk\"ahler embeddings are completely geodesic}

In assumptions of \ref{_NX_hyperholo_Proposition_}, 
consider the bundle $TM\restrict X$ with metrics and connection
induced from $TM$. There is a natural exact sequence of 
holomorphic vector bundles over $(X, L)$:
\[
  0\arrow TX \arrow TM\restrict X \arrow NX \arrow 0 .
\]
The natural connection in each of these bundles is hyperholomorphic
(\ref{_NX_hyperholo_Proposition_}).

\hfill

\proposition \label{_NX_splits_for_hype_Proposition_}
Let $X$ be a hyperk\"ahler manifold, not necessary compact, $L$
induced complex structure, and 
\begin{equation}\label{_holo_seque_Equation_}
    0\arrow E_1\arrow E_2 \arrow E_3\arrow 0
\end{equation}
 be an exact sequence
of holomorphic vector bundles. Let $g_2$ be a Hermitian structure 
on $E_2$. Consider the metrics $g_1$, $g_2$, $g_3$ and the connections
$\nabla_1$, $\nabla_2$, $\nabla_3$ on $E_1$, $E_2$, $E_3$ 
induced by $g_2$. Assume that $\nabla_1$, $\nabla_2$, $\nabla_3$
are hyperholomorphic. Then the exact sequence
\eqref{_holo_seque_Equation_} splits, and moreover,
the orthogonal decomposition $E_2 = E_1 \oplus E_1^\bot$ is preserved by
 the connection $\nabla_2$.

\hfill

{\bf Proof:} The same statement is well known 
for the Yang--Mills connections: every exact sequence of 
holomorphic bundles with compatible Yang--Mills metrics
splits (see Appendix to this section). 
By \ref{_hyperholo_Yang--Mills_Proposition_}, 
hyperholomorphic connections
are always Yang--Mills. \endproof

\hfill

\definition
Let $N \stackrel i \hookrightarrow M$ be an embedding of
hyperk\"ahler manifolds. We say that $i$ is a {\bf hyperk\"ahler
embedding} if $i$ is compatible with the quaternionic structure
and Riemannian metric.

\hfill

\corollary\label{_hype_embe_comple_geode:Corollary_}
Let $N \stackrel i \hookrightarrow M$ be a hyperk\"ahler embedding.
Then $i$ is completely geodesic.

{\bf Proof:} Follows directly from 
\ref{_comple_geodesi_basi_Proposition_} and 
\ref{_NX_splits_for_hype_Proposition_} \endproof

\hfill

\subsection{Appendix: every exact 
sequence of Yang--Mills bundles splits.}

\hspace{6mm}
\proposition \label{_YM_exact_split_Proposition_}
Let $X$ be a K\"ahler manifold, not necessary compact, and
\begin{equation}\label{_holo_seque2_Equation_}
    0\arrow E_1\arrow E_2 \arrow E_3\arrow 0
\end{equation}
 be an exact sequence
of holomorphic vector bundles. Let $g_2$ be a Hermitian 
structure on $E_2$.
Consider the induced metrics $g_1$, $g_2$, $g_3$. 
Assume that either $g_1$ or $g_3$ is Yang--Mills. 
Then the exact sequence
\eqref{_holo_seque_Equation_} splits, and moreover,
the orthogonal complement $E_1^\bot\subset E_2$ is preserved by
 the connection $\nabla_2$.

\hfill

{\bf Proof:} Consider the second fundamental form 
\[ A\in \Lambda^{0,1}(X, \Hom(E_1,E_3)) \]
of the exact sequence \eqref{_holo_seque2_Equation_}.
The curvatures $\Theta_i$ of $E_i$ are expressed through
$A$ as follows (\cite{_Griffi_Harri_}): 
\begin{equation}\label{_Theta_1_through_seco_Equation_}
   \Theta_1 = \Theta_2 \restrict{E_1} + {}^t\bar A\wedge A,
\end{equation}
\begin{equation}\label{_Theta_3_through_seco_Equation_}
   \Theta_3 = \Theta_2 \restrict{E_3} - A\wedge{}^t\bar A.
\end{equation}
The 2-forms ${}^t\bar A\wedge A$, $A\wedge{}^t\bar A$ are {\bf positive}
unless $A=0$. Thus, the endomorphisms
$\Lambda \left(A\wedge{}^t\bar A\right)\in End(E_1)$, 
$\Lambda \left({}^t\bar A\wedge A\right)\in End(E_3)$ have positive trace
(again, unless $A=0$). By our assumption,
$\Lambda (\Theta_2)=0$ and either $\Lambda(\Theta_1)=0$
or $\Lambda(\Theta_3)=0$. Applying $\Lambda$ to both sides
of \eqref{_Theta_1_through_seco_Equation_} and
\eqref{_Theta_3_through_seco_Equation_}, we obtain that either
$\Lambda \left(A\wedge{}^t\bar A\right)=0$ or
$\Lambda \left({}^t\bar A\wedge A\right)=0$. 
Therefore, $A=0$ and the exact sequence
\eqref{_holo_seque2_Equation_} splits.
\endproof


\section[Triholomorphic sections of hyperholomorphic bundles
and deformations of trianalytic submanifolds.]
{Triholomorphic sections \\of hyperholomorphic bundles\\
and deformations of trianalytic submanifolds.}
\label{_triholo_Section_}


In the previous section, we proved that hyperk\"ahler embeddings are
completely geodesic. In this section, we show that, furthermore,
the deformational results of Section
\ref{_comple_geode_defo_Section_} are fully applicable
to the deformations of trianalytic submanifolds.

\subsection{Triholomorphic sections of normal bundle}

\hfill

\definition 
Let $M$ be a hyperk\"ahler manifold, $B$ a vector bundle
equipped with a hyperholomorphic connection, and $\alpha$ a section
of $B$. Then $\alpha$ is called {\bf triholomorphic} 
if for each induced complex structure $L$ on $M$, $\alpha$ 
is a holomorphic section of $B$ considered as a holomorphic
bundle over $(M, L)$. 

\hfill

Let $M$ be a hyperk\"ahler manifold, not necessarily compact, and
$X\subset M$ a trianalytic submanifold, not necessarily closed.
Fix an induced complex structure $I$ on $M$. Let 
\[ 
   \left( \pi:\; \c X \arrow S, 
   \phi:\; \c X \arrow M, \phi(\pi^{-1}(s_0)) = X \right)
 \]
be a family of deformations of $(X,I) \subset (M, I)$
(see Section \ref{_comple_geode_defo_Section_} for details). 
Let $\gamma:\; [0,1] \arrow S$ be a real analytic path in $S$, such
that $\gamma(0) = s_0$. Assume that for all $t\in [0,1]$, the
submanifold $X_t = \phi(\gamma^{-1}(t))\subset M$ is trianalytic.
Consider the normal bundle $NX$ to $X$, with the metric and
connection induced from $M$. By \ref{_NX_hyperholo_Proposition_}, 
$NX$ is hyperholomorphic. Let $\eta$ be the section of 
$NX$ corresponding to $\frac{d\gamma}{dt}$ as in 
Subsection \ref{_norma_vecto_Subsection_}.

\hfill

\proposition \label{_norma_triholo_Proposition_}
In the above assumptions, $\eta$ is a triholomorphic section of $NX$.

\hfill

{\bf Proof:} Let $L$ be an induced complex structure, and $(M, L)$,
$(X, L)$ the manifolds $M$ and $X$ considered as complex 
manifolds with the complex structure $L$. The normal bundle
$N(X,L)$ is naturally identified with $NX$ as a real vector 
bundle. Therefore $\eta$ can be considered as a section
of $N(X,L)$. By \ref{_real_ana_indu_on_hype_equiva_Proposition_}
(see Appendix to this section), $\c X$ can be considered
as a real analytic deformation of $(X, L)$.
{}From \ref{_norma_sec_holom_Lemma_} 
it is clear that $\eta$ is holomorphic
as a section of $N(X,L)$. This proves  
\ref{_norma_triholo_Proposition_}. \endproof

\subsection{Triholomorphic sections are parallel}

\hfill

\proposition\label{_triholo_parallel_Proposition_}
Let $M$ be a hyperk\"ahler manifold, not necessary compact,
and $B$ a vector bundle with a hyperholomorphic connection $\nabla$.
Let $\nu$ be a trianalytic section of $B$. Then, $\nu$ is parallel:
\begin{equation}\label{_nu_para_Equation_}
\nabla\nu=0
\end{equation}

{\bf Proof:} Let $L$ be an induced complex structure. Since 
$\nu$ is triholomorphic, $\bar \6_L \nu =0$, where
$\bar \6_L:\; B \arrow B\times \Lambda^{0,1}_L(M)$ is the
$(0,1)$-part of the connection, taken with respect to $L$.
Taking $L=I, -I$, we obtain 
\begin{equation} \label{_bar6_I+bar6_-I_Equation_}
\bar \6_I + \bar\6_{-I}(\nu)=0.
\end{equation}
On the other hand, $\bar\6_{-I}= \6_I$, where $\6_I$ is the 
$(1,0)$-part of $\nabla$ taken with respect to $I$. Thus,
\[ \bar \6_I + \bar\6_{-I} =\bar \6_I + \6_{I} =\nabla. \]
{}From \eqref{_bar6_I+bar6_-I_Equation_} 
we obtain that $\nabla(\nu)=0$. \endproof

\hfill

{}From \ref{_norma_triholo_Proposition_} 
and \ref{_triholo_parallel_Proposition_}, 
we obtain that the section $\eta \in NX$ 
is parallel with respect to the connection. Then, 
\ref{_conne_in_fami_of_comple_geo_Proposition_} (i)
can be applied to the following effect.

\hfill

\corollary \label{_norma_sec_para_for_hype_from_it_Corollary_}
Let $M$ be a hyperk\"ahler manifold, not necessarily compact, and
$X\subset M$ a trianalytic submanifold, not necessarily closed,
and
\[ 
   \left( \pi:\; \c X \arrow S, 
   \phi:\; \c X \arrow M, \phi(\pi^{-1}(s_0)) = X \right)
 \]
be a real analytic family of deformations of $X\subset M$
(see Section \ref{_comple_geode_defo_Section_} for details). 
Let $\gamma:\; [0,1] \arrow S$ be a real analytic path in $S$, such
that $\gamma(0) = s_0$. Assume that for all $t\in [0,1]$, the
submanifolds $X_t = \phi(\gamma^{-1}(t))\subset M$ are trianalytic.
Fix an induced complex structure $I$ on $M$.
Let $U_t\subset (X_t, I)$, $t\in [0,1]$ be the subsets constructed 
in Subsection \ref{_conne_Subsection_}\footnote{By 
\ref{_hype_embe_comple_geode:Corollary_},
trianalytic submanifolds are completely geodesic, and thus,
\ref{_conne_in_fam_of_comple_geode_Proposition_} can be applied.}
with the corresponding holomorphic isomorphisms 
$\Psi_t:\; U_0\arrow U_t$. Then, assumptions of
\ref{_conne_in_fami_of_comple_geo_Proposition_} (i) hold.
Thus,
for all $t, t'\in [0,1]$,
the maps $\Psi^t_{t'}$ are isometries. Moreover,
$\Psi^t_{t'}$ are compatible with the hyperk\"ahler structure.

\hfill

{\bf Proof:} 
Assumptions of \ref{_conne_in_fami_of_comple_geo_Proposition_} (i)
hold by \ref{_norma_triholo_Proposition_} 
and \ref{_triholo_parallel_Proposition_}.
The maps $\Psi^t_{t'}$ are isometries by 
\ref{_conne_in_fami_of_comple_geo_Proposition_} (i).
The maps $\Psi^t_{t'}$ are compatible with the 
hyperk\"ahler structure because the they are obtained by integrating
a certain connection in the family $\c X$. This connection
is constructed from the Hermitian metric, and thus, does
not depend from the choice of induced complex structure. Taking
different induced complex structures,
we obtain the same maps $\Psi^t_{t'}$. 
Thus, $\Psi^t_{t'}$ is holomorphic
with respect to each of induced complex structures
(see also \ref{_norma_sec_holom_Lemma_}).
\endproof

\subsection{Appendix: real analytic structures on hyperk\"ahler
manifolds}

Consider the real analytic structures on a given hyperk\"ahler
manifold arising from the different induced complex structures.
We prove that these real analytic structures are equivalent.

\hfill

\proposition \label{_real_ana_indu_on_hype_equiva_Proposition_} 
Let $M$ be a hyperk\"ahler manifold, $I_1$, $I_2$ induced complex
structures. Let $(M, I_1)$ and $(M, I_2)$ be the corresponding
complex manifolds, and $(M, I_1)_\R$, $(M, I_2)_\R$ be the real
analytic manifolds underlying $(M, I_1)$, $(M, I_2)$. Consider the
tautological map $(M, I_1)_\R\stackrel \phi \arrow (M, I_2)_\R$.
Then $\phi$ is compatible with the real analytic structure.

\hfill

{\bf Proof:} Consider the {\bf twistor space} for $M$ (see
\cite{_Besse:Einst_Manifo_}), $\Tw(M)$, with the natural
holomorphic map $\pi:\; \Tw(M) \arrow \C P^1$. Let $Sec(M)$
be the space of holomorphic sections of the map $\pi$.
Then $Sec(M)$ is identified naturally with an open subspace
of a Douady space for $\Tw(M)$, and thus, has a natural 
complex structure. According to \cite{_HKLR_} (see \cite{_NHYM_} 
for details), $Sec(M)$ is a complexification of $(M,I)$, in the
sense of Grauert. In other words, the complex valued real analytic functions
on $(M,I)$ are naturally identified with the germs of complex analytic
functions on $Sec(M)$. Since $Sec(M)$ is defined independently
from the choice of an induced complex structure, the 
tautological map $\phi$ is an equivalence.
This proves \ref{_real_ana_indu_on_hype_equiva_Proposition_}.
 \endproof

\hfill

\ref{_real_ana_indu_on_hype_equiva_Proposition_} implies that we
may speak of a real analytic manifold underlying a given hyperk\"ahler
manifold. 

\hfill

\corollary \label{_Riema_on_hype_real_ana_Corollary_}
Let $M$ be a hyperk\"ahler manifold. Consider the Riemannian
form $g$ on $M$ as a section of the real analytic bundle of 
symmetric 2-forms. Then $g$ is real analytic.

\hfill

{\bf Proof:} Let $I, J, K$ be the induced complex structures
which form the standard basis in quaternions, and 
$\omega_I$, $\omega_J$, $\omega_K$ be the corresponding 
K\"ahler forms. Then $\Omega:=\omega_J + \1\omega_K$ is the natural
holomorphically symplectic form on $(M,I)$, and as such,
$\Omega$ is real analytic. Then, its real part $\omega_J$ is
also real analytic. Since the complex structure operator $J$
is real analytic, we obtain that the form
$g(\cdot,\cdot):= - \omega(\cdot, J\cdot)$ is also real analytic.
\endproof


\section{Douady spaces for trianalytic cycles and real analytic structure}
\label{_Douady_Section_}


\subsection{Real analytic structure on the Douady space.} 

In this section, we consider the Douady space $D_M(X)$
for a subvariety $X$
of a compact complex manifold $M$ equipped with
 a hyperk\"ahler structure.
We prove that, when $X$ is trianalytic, $D_M(X)$ is {\bf
hypercomplex} (\ref{_hypercomplex_Definition_}). 
We also show that the real analytic
variety underlying $D_M(X)$ does not change if we replace
a complex structure on $M$ by another induced complex structure.
These results are technical and we use them mainly to 
simplify the exposition. 

\hfill

Let $M$ be a compact hyperk\"ahler manifold, and $X\subset M$ a 
closed trianalytic subvariety. For each induced complex
structure $L$, consider $(X, L)$ as a complex subvariety of $(M, L)$.
Let $D_L(X)$ be the Douady deformation space of $(X, L)$ in
$(M, L)$. The points $[X']$ of $D_L(X)$ correspond to the
subvarieties $X'\subset (M, L)$ which are deformation of $X$.
By \ref{_G_M_invariant_implies_trianalytic_Theorem_}  (see also
\cite{Verbitsky:Symplectic_II_}), every such $X'$ is trianalytic. 

Let $I_1$, $I_2$ be induced complex structures on $M$.
The following proposition shows that a trianalytic subvariety
$X'$ can be obtained as a deformation of $(X, I_1)$ in $(M, I_1)$
if and only if $X'$ can be obtained
as a deformation of $(X, I_2)$ in $(M, I_2)$.

\hfill

\proposition \label{_X_in_Douady_indep_from_I_1_Proposition_} 
Let $I_1$, $I_2$ be induced complex
structures on a compact hyperk\"ahler manifold
$M$, and $X\subset M$ a closed trianalytic subvariety. 
Consider the corresponding Douady spaces $D_{I_1}(X)$,
$D_{I_2}(X)$. Then, for a trianalytic subvariety $X'\subset M$,
$[X'] \in D_{I_1}(X)$ if and only if $[X']\in D_{I_2}(X)$.

\hfill

{\bf Proof:} Let us recall the notion of {\bf degree} of
a trianalytic subvariety.
Let $M_1$ be a K\"ahler manifold.
By {\bf degree} of a subvariety $X \subset M_1$ we understand a number
\[ \deg X := \int_X \omega^{\dim_\C X}, \]
where $\omega$ is the K\"ahler form. If $X\subset M$ is 
a trianalytic subvariety of a hyperk\"ahler manifold, we can
associate a number $\deg X$ to each of induced complex structures.
In \cite{Verbitsky:Symplectic_I_}, we prove that $\deg X$
is in fact independent from the choice of induced complex structure.
This enables us to speak of {\bf degree} of a closed trianalytic
subvariety of a compact complex manifold.

We return to the proof of
\ref{_X_in_Douady_indep_from_I_1_Proposition_}.
Let $D_1$ be the union of all components of the 
Douady space for $(M,I_1)$. By \cite{_Lieberman_},
\cite{_Fukjiki_Kahler_}, $D_1$ is compact;
in particular, $D_1$ has a finite number of connected components.
Let $D_1 = \cup D^i_1$, $i\in \Upsilon$ 
be the decomposition of $D$ unto a union of its
connected components. 

Let $X'\subset M$ be a closed trianalytic subvariety of $M$
such that $[X'] \in D_{I_2}(X)$. There exists a  real analytic
path $\gamma:\; [0,1] \arrow D_{I_2}(X)$ joining $[X]$ and $[X']$.
For all $t\in [0,1]$, the point $\gamma(t)$ lies in $D_1^i$ for some
$i\in \Upsilon$. Let $S^i\subset [0,1]$ be the set of all
$t\in [0,1]$ such that $\gamma(t)\in D_1^i$. We are going to
show that $S^i$ are compact for all $i$. This will clearly imply
that for all $i$ except one, $S^i$ is empty, thus proving that
$[X]$ and $[X']$ lie in the same component of $D_1$.

Let $D_\R$ be the set of all real analytic deformations of $X$ in $M$
with natural topology. The forgetful map $\psi:\; D_1 \arrow D_\R$
is continous and injective. Since $D_1$ is compact, $\psi$
is a closed embedding. The composition 
$\gamma\circ \psi:\; [0,1] \arrow D_\R$ is obviously continous,
and thus, $\gamma:\; [0,1] \arrow D_1$ is also continous. 
This implies that all $S_i$, $i\in \Upsilon$ are 
closed, thus proving that all $S_i$ except one are empty.
\ref{_X_in_Douady_indep_from_I_1_Proposition_} is proven.
\endproof

\hfill

Let $D_{I_1}(X)_\R$, $D_{I_2}(X)_\R$ be the real analytic 
varieties underlying $D_{I_1}(X)$, $D_{I_2}(X)$. 
\ref{_X_in_Douady_indep_from_I_1_Proposition_} gives a 
tautological bijection \[ \psi:\; D_{I_1}(X)_\R\arrow D_{I_2}(X)_\R.\]

\hfill

\proposition\label{_Douady_real_ana_inde_from_indu_Proposition_}
The map $\psi:\; D_{I_1}(X)_\R\arrow D_{I_2}(X)_\R$
is an isomorphism of real analytic varieties.

\hfill

{\bf Proof:} Let $I$ be an induced complex structure on $M$.
For a real analytic function $f$ on $M$, consider the function
$\hat f:\; D_I(X) \arrow \R$,
\[ [X'] \stackrel{\hat f}\arrow \int_{X'} f \cdot \Vol X', \]
with $\Vol X'$ the volume form on $X'$. For
all $f$, the function $\hat f$ is real analytic. From the definition of
Douady spaces it might be seen that converging 
power series of different $\hat f$ generate the sheaf
$\c A_I$ of real analytic functions on $D_I(X)$. For an induced
complex structure $I'$, the same set $\{ \hat f\}$
generates the sheaf $\c A_{I'}$ 
of real analytic functions on $D_{I'}(X)$
(\ref{_real_ana_indu_on_hype_equiva_Proposition_}). Thus,
the sheaves $\c A_I$ and $\c A_{I'}$ coinside.
\endproof

\hfill

\ref{_Douady_real_ana_inde_from_indu_Proposition_}
shows that we may speak of the real analytic variety underlying
$D_M(X)$ without specifying an induced complex structure $I$.

\hfill

\definition \label{_hypercomplex_Definition_}
Let $Y$ be a real analytic variety equipped with
three complex structures $I$, $J$ and $K$. Assume that
for every point $y\in Y$, the action of $I$, $J$, $K$ on the
Zariski tangent space $T_yY$ satisfies $I\circ J = - J\circ I =K$.
Then $Y$ is called {\bf a hypercomplex variety}.

\hfill

\remark \label{_Douady_hyperc_Remark_}
\ref{_Douady_real_ana_inde_from_indu_Proposition_}
immediately implies that the induced complex structures on $M$ 
equip the Douady space $D_M(X)$ with a hypercomplex structure.


\subsection{Appendix: isometric embeddings are hyperk\"ahler.}


Let $M$ be a compact hyperk\"ahler manifold, $X\subset M$
a closed trianalytic subvariety.
There is an alternative way to describe the Douady space
$D(X)$. The \ref{_isome_embe_Proposition_} 
below describes $D(X)$ in terms of of the space of isometric embeddings
$\nu:\; X \arrow M$. 

\hfill

Let $M$ be a compact K\"ahler manifold, and $X\subset M$
be a closed analytic subvariety. One can make sense of the
{\bf fundamental class} of $X$, $[X] \in H^{\dim_\R(X)}(M)$,
lying in homology of $M$.
Consider $X$ with an induced structure of a metric space,
and let $X \stackrel \nu \arrow M'$ be an isometric embedding,
where $M'$ is a compact Riemannian manifold. The same argument which
allows us to define the fundamental class of $X$, allows
us to define the fundamental class of $\nu(X)$. 

\hfill

\proposition \label{_isome_embe_Proposition_}
Let $M$ be a closed analytic subvariety, $X\subset M$ a closed 
trianalytic subvariety, and $\nu:\; X \arrow M$ an isometric embedding.
Assume that the fundamental class of $\nu(X)$ is equal to the fundamental
class of $X$. Then $\nu(X)$ is a trianalytic subvariety of $M$.

\hfill

{\bf Proof:} We use notation 
introduced in Subsection \ref{_SU(2)-inv=>triana_Subsection_}.
For each induced complex structure $I$,
we have 
\begin{equation}\label{_deg_same_appl_nu_Equation_}
   \deg_I (X) - \deg_I (\nu(X)),
\end{equation} 
because $\nu(X)$ has the same fundamental class as $X$.
Since $X$ is trianalytic,
applying Wirtinger's inequality (\ref{_Wirti_for_Kahle_Theorem_}), 
we obtain that
\begin{equation} 
   2^n \int_{X_{ns}} \Vol X_{ns} = \int_{X_{ns}} \omega^n,
\end{equation}
where $\omega$ is the K\"ahler class of $(M,I)$, $n=\frac{1}{2}\dim_\R(X)$
and $X_{ns}$ is the nonsingular part of $(X,I)$.
By definition, we have
\begin{equation}
   \int_{X_{ns}} \omega^n = \deg_I (X), \ \ 
   \int_{\nu(X_{ns})} \omega^n = \deg_I (\nu(X)).
\end{equation}
Since $\nu$
is an isometry, and $\Vol(X)$ is an invariant of a metric, 
we have 
\begin{equation} \label{_Vol_same_appl_nu_Equation_}
   \int_{X_{ns}} \Vol X_{ns} = 
   \int_{\nu(X_{ns})} \Vol \left(\nu(X_{ns}) \right).
\end{equation}
Together the equations 
\eqref{_deg_same_appl_nu_Equation_}--\eqref{_Vol_same_appl_nu_Equation_}
give 
\[ 
   2^n\int_{X_{ns}} \Vol X_{ns} = \int_{\nu(X_{ns})} \omega^n. 
\]
Applying Wirtinger's inequality (\ref{_Wirti_for_Kahle_Theorem_}) once
again, we obtain that $\nu(X)$ is complex analytic with respect to $I$.
This proves \ref{_isome_embe_Proposition_}. \endproof

\hfill

We just proved \ref{_iso_intro:Theorem_} (iii). 
As another application of \ref{_isome_embe_Proposition_},
we give a direct proof that $D(X)$ is compact.

\hfill  

\corollary
Let $M$ be a compact hyperk\"ahler manifold, 
$X\stackrel {\bar i}\hookrightarrow  M$ be a closed trianalytic subvariety
and $D(X)$ its Douady space. Then $D(X)$ is compact.

\hfill

{\bf Proof:} Let $X_{ns}$ be the non-singular part of $X$,
and $i:\; X_{ns} \hookrightarrow M$ the natural embedding.
Then $i$ is an isometry. Consider a deformation $\nu$ of $i$ in class of
isometries. The arument proving \ref{_isome_embe_Proposition_}
shows that the closure of $\nu(X_{ns})$ is trianalytic in $M$.
Consider the space $\underline{Is(X_{ns})}$ of isometries 
$X_{ns}\stackrel \nu\hookrightarrow M$. Let $Is(X_{ns})$ be
a connected component of this space, containing $i$. There
is a natural continuous surjection 
$p:\; Is(X_{ns}) \arrow D(X)$, which maps $\nu \in Is(X_{ns})$
to a closure of $\nu(X_{ns})$ in $M$. To prove that $D(X)$ is
compact it suffices to show that $Is(X_{ns})$ is compact.
This is an implication of the following general statement, which is clear.

\hfill

\claim 
Let $X$, $M$ be Riemannian manifolds, $M$ compact. Let
$Is(X, M)$ be the space of isometries (maps which preserve the 
geodesic distance) from $X$ to $M$.  Then $Is(X, M)$ is compact.

\endproof


\section[Connections in the families of trianalytic subvarieties]
{Connections in the families of trianalytic \\subvarieties}
\label{_Conne_in_fami_Section_}


\subsection{Introduction}
\label{_intro_conne_Subsection_}

Let $M$ be a compact hyperk\"ahler manifold, and
\begin{equation} \label{_fami_Equation_}
   \left( \pi:\; {\c X} \arrow S, \phi:\; {\c X} \arrow M 
   \right )
\end{equation}
a real analytic family of subvarieties, not necessarily closed. 
Assume that for all $s\in S$, the fiber $X_s = \phi(\pi^{-1} (s))$
is trianalytic in $M$.
By \ref{_conne_in_fam_of_comple_geode_Proposition_} 
and \ref{_hype_embe_comple_geode:Corollary_}, the family 
${\c X}$ is then equipped with a connection
$\nabla:\; N_\pi X \arrow T {\c X}$. It is natural 
to ask whether this connection is {\bf flat} 
(see \ref{_flat_conne_in_fami_Definition_}). 
The answer is affirmative, 
under certain additional assumptions.

\hfill

\definition \label{_admits_co_Definition_}
We say that the family \eqref{_fami_Equation_} {\bf admits 
a compactification} if the following conditions hold.

\begin{description}
\item[(i)]
For each $s \in S$, the closure $\bar X_s$ of the fiber 
$X_s$ in $M$ is a trianalytic subvariety of $M$.

\item [(ii)] For all $s\in S$, $\bar X_s$ lie in the same 
component of the Douady space.

\end{description}

\hfill

We show that in assumptions of 
\ref{_admits_co_Definition_}, the connection
\[ \nabla:\; N_\pi X \arrow T {\c X}\] is indeed flat
(\ref{_conne_in_triana_flat_Theorem_}). 
However, we don't know whether it is flat
for general families of trianalytic submanifolds.

\hfill


\subsection{Flat connections and curvature.}


\hfill

\definition
(curvature of a connection) 
Let $\pi: {\c X} \arrow S$ be a family of manifolds,
$T_\pi X \subset T{\c X}$ be the fiberwise tangent bundle and 
$N_\pi X = T{\c X} / T_\pi X$ be a fiberwise normal bundle. 
Let $\nabla:\; N_\pi X \arrow T{\c X}$ be a connection in ${\c X}$.
Then the {\bf curvature} of $\nabla$ is the following
tensor 
$\Theta\in \Hom_\R\left ( \Lambda^2_\R \left( N_\pi X\right), 
T_\pi{\c X}\right)$.
For two sections $a,b \in N_\pi X$, consider the corresponding
vector fields $\nabla a$, $\nabla b\in T{\c X}$. Consider the projection
$T {\c X} \stackrel {pr}{\arrow} T{\c X} / \nabla N_\pi X = T_\pi X$.
Let $\Theta(a,b):= pr([\nabla a, \nabla b])$.
Clearly, $\Theta$ is a tensor. 

\hfill

\proposition
Let $\left( \pi:\; {\c X} \arrow S \right)$ 
be a family of manifolds equipped with a connection $\nabla$.
Then the following conditions are equivalent.

\begin{description}

\item[(i)] The connection $\nabla$ is flat.

\item[(ii)] Let $s_1, s_2\in S$ be a pair of points,
and $\gamma, \gamma':\; [0,1] \arrow S$ real analytic paths.
Let $U_0 \subset X_{s_0}$, $U_1 \subset X_{s_1}$ be the open
subsets such that the connection $\nabla$ might be integrated to a map
$\Psi:\; U_{s_0} \arrow U_{s_1}$ along $\gamma$, and
$U_0' \subset X_{s_0}$, $U_1' \subset X_{s_1}$ be 
such that $\nabla$ might be integrated to a map
$\Psi':\; U_{s_0}' \arrow U_{s_1}'$ along $\gamma'$. Then
$\Psi$ coinsides with $\Psi'$ in the intersection
$U_{s_0} \cap U'_{s_0}$.

\end{description}

{\bf Proof:} Well known. \endproof

\hfill

\definition \label{_flat_conne_in_fami_Definition_} 
 Let ${\c X} \arrow S$ be a family of manifolds, equipped with a
connection $\nabla$. Then $\nabla$ is called {\bf flat}
if its curvature is zero.


\subsection{Curvature of holomorphic connections}


Let $\pi:\; {\c X} \arrow S$ be a complex analytic family
of manifolds. Then the vector bundles $T{\c X}$ and
$N_\pi X$ are equipped with a natural holomorphic structure.
The connection $\nabla:\; N_\pi X \arrow T {\c X}$ is called
{\bf holomorphic} if the map $\nabla:\; N_\pi X \arrow T{\c X}$
is holomorphic. 

\hfill

\claim \label{_holo_conne_curva_C-line_Claim_}
Let $\pi:\; {\c X} \arrow S$ be a family of 
manifolds equipped with with a holomorphic connection
$\nabla$. Consider the curvature 
\[ \Theta\in \Hom_\R\left 
   ( \Lambda^2_\R \left(N_\pi X\right), T_\pi{\c X}\right).
\]
Then $\Theta$ is $\C$-linear; in other words, $\Theta$ belongs
to the subspace

\[ \Hom_\C\left ( \Lambda^2_\C \left(N_\pi X\right), T_\pi{\c X}\right)
   \subset \Hom_\R\left 
    ( \Lambda^2_\R \left(N_\pi X\right), T_\pi{\c X}\right).
\]

{\bf Proof:} Clear. \endproof


\subsection{Flat connections in a family of trianalytic varieties. }


Let $M$ be a compact hyperk\"ahler manifold, $X\subset M$ 
a closed trianalytic subvariety, $D(X)$ its Douady space
and 
\[ \left( \underline \pi: \; \underline{{\c X}} \arrow D(X), 
   \underline \phi:\; \underline{{\c X}} \arrow M \right )
\]
a universal family of subvarieties corresponding to $D(X)$.
Let ${\c X}$ be the union of all smooth points in all fibers
of $\underline \pi$, and $\pi$, $\phi$ be the restrictions of
$\underline \pi$ and $\underline \phi$ to ${\c X}$.
Consider the corresponding family of manifolds
\[ \left(  \pi: {{\c X}} \arrow D(X), 
    \phi {{\c X}} \arrow M \right ).
\]
We are in the same situation as described in 
Subsection \ref{_intro_conne_Subsection_}. 
Thus, the family ${\c X}$ is equipped with a natural 
connection $\nabla$.

\hfill

\theorem \label{_conne_in_triana_flat_Theorem_}
In these assumptions, $\nabla$ is flat.

\hfill

{\bf Proof:} The proof of \ref{_conne_in_triana_flat_Theorem_}
takes the rest of this section. 

\hfill

Clearly, $N_\pi X = \pi^* TS$. By \ref{_Douady_hyperc_Remark_}, 
$TS$ has a natural quaternionic action, and thus, an action of
$SU(2) \subset {\Bbb H}^*$. The bundle $T_\pi X$ is also equipped 
with a natural action of quaternions. This endows the bundle
$\Hom_\R\left ( \Lambda^2_\R \right(N_\pi X\left), T_\pi{\c X}\right)$
with an action of $SU(2)$.

\hfill

\lemma  \label{_curva_SU(2)_inva_Lemma_}
Let $\Theta$ be a curvature of $\nabla$, 
$\Theta \in\Hom_\R
\left ( \Lambda^2_\R \left( N_\pi X\right), T_\pi{\c X}\right)$.
Then $\Theta$ is $SU(2)$-invariant.

\hfill

{\bf Proof:} Let $U(1) \stackrel {\tilde i}\hookrightarrow SU(2)$ be 
an embedding corresponding to an algebra embedding 
$\C \stackrel i\hookrightarrow {\Bbb H}$. Let $I$ be an
induced complex structure associated with $i$. The
connection $\nabla$ is holomorphic with respect to $I$.
Therefore $\Theta$ is $\C$-linear
with respect to the action of $\C$ given by 
$\C \stackrel i \hookrightarrow {\Bbb H}$
(\ref{_holo_conne_curva_C-line_Claim_}). 
Thus, $\Theta$ is $U(1)$-invariant with respect
to the action of $U(1)$ given by 
$\tilde i:\; U(1)\hookrightarrow SU(2)$.
The group $SU(2)$ is generated by the images of $\tilde i$ for all
algebra embeddings $i:\; \C \hookrightarrow {\Bbb H}$. Thus,
$\Theta$ is $SU(2)$-invariant. \endproof

\hfill

Consider $N_\pi X$, $T_\pi X$ as representation of $SU(2)$. These representations
are of weight 1. The bundle 
\[ 
   \Hom_\R\left ( \Lambda^2_\R \left(N_\pi X\right), T_\pi{\c X}\right)
\] 
is 
again a representation of $SU(2)$. Since 
$N_\pi X$, $T_\pi X$ are of weight one, 
$\Hom_\R\left ( \Lambda^2_\R \left(N_\pi X\right), T_\pi{\c X}\right)$ is 
a direct sum of representations of weight 3 and 1. Thus,
this bundle cannot have non-zero $SU(2)$-invariant sections.
Therefore, ($\Theta$ is $SU(2)$-invariant) 
implies ($\Theta =0$).
\ref{_conne_in_triana_flat_Theorem_} is proven. \endproof


\section[Isometries of trianalytic subvarieties 
of a compact hy\-per\-k\"ah\-ler manifold]{Isometries 
of trianalytic subvarieties \\ 
of a compact hyperk\"ahler manifold}
\label{_isome=>holo_Section_}


\subsection{Premises}

Let $M$ be a compact hyperk\"ahler manifold, $X\subset M$ a 
closed trianalytic subvariety.
{}From \ref{_conne_in_fami_of_comple_geo_Proposition_},
\ref{_conne_in_triana_flat_Theorem_} and 
\ref{_norma_sec_para_for_hype_from_it_Corollary_},
we immediately obtain the following theorem.

\hfill

\theorem \label{_triana_subse_isome_Theorem_}
Let $M$ be a compact hyperk\"ahler manifold, $X\subset M$
a closed trianalytic subvariety, $D(X)$ the corresponding 
Douady space (\ref{_Douady_hyperc_Remark_}). 
Let $\gamma:\; [0,1] \arrow D(X)$ 
be a real analytic path. Denote by $X_t$ 
the subvarieties corresponding to points $\gamma(t)\in D(X)$.
Consider $X_t$ as metric spaces with the metric induced
from $M$.
Then, for each $t_1$, $t_2\in [0,1]$,
there exist a natural isometry
$\Psi^{t_1}_{t_2}:\; X_{t_1} \arrow X_{t_2}$, mapping
non-singular points to non-singular points, and
acting compatible with the hyperk\"ahler structure 
on $X_{t_1}^{ns}$.\footnote{The non-singular part $X_t^{ns}$ 
is naturally equipped with a structure of a hyperk\"ahler manifold.}
This isometry depends only on the homotopy class of $\gamma$.

\endproof

\hfill

This almost finishes the proof \ref{_iso_intro:Theorem_} (i).
It remains to prove the following
theorem.

\hfill

\theorem \label{_triana_subse_comple_ana_Theorem_}
In assumptions of \ref{_triana_subse_isome_Theorem_},
let $I$ be an induced complex structure on $M$. Then
the isometry $\Psi^{t_1}_{t_2}:\; X_{t_1} \arrow X_{t_2}$
is compatible with the complex analytic structure induced
by $I$.

\hfill

This section is taken fully with the proof of 
\ref{_triana_subse_comple_ana_Theorem_}


\subsection{Homeomorphisms of complex varieties and normalization.}


\hfill

\claim \label{_conti_fu_extend_to_holo_Claim_}
(\cite{_Grauert_}) Let $X$ be a normal complex analytic variety, $U$
a dense open subset in $X$ and $f:\; U \arrow \C$ a bounded 
holomorphic function. Then $f$ can be extended to a holomorphic function
on $X$.

\endproof

\hfill

This statement has an immediate corollary.

\hfill

\corollary \label{_conti_map_norma_holo_Corollary_}
Let $\phi:\; X \arrow Y$ be a continous map of complex analytic
varieties. Assume that $\phi$ is holomorphic in an open
dense subset $U\subset X$. Then, if $X$ is normal, $\phi$
is holomorphic.

\hfill

{\bf Proof:} Let $x\in X$ be an aritrary point, $V$ its neighbourhood,
sufficiently small. Taking 
coordinates in a neighbourhood $W$ of $\phi(V)$ and applying
\ref{_conti_fu_extend_to_holo_Claim_}, we find that
$\phi\restrict {U\cap V}$ extends continously to a 
holomorphic map $\tilde \phi$ from $U\cap V$ to $W$. 
Since $U\cap V$ is dense in $V$, $\tilde \phi$
coinsides with $\phi\restrict V$. Thus, $\phi$ is
holomorphic. \endproof

\hfill

Thus, in the situation of \ref{_triana_subse_comple_ana_Theorem_},
were $X_t$ normal, $\Psi^{t}_{t'}$ would have been holomorphic
and \ref{_triana_subse_comple_ana_Theorem_} would have been proven.
Unfortunately, we have no means to show that $X_t$ is normal.
However, from \ref{_conti_map_norma_holo_Corollary_}
we obtain some information about maps of arbitrary 
varieties too.

\hfill

\corollary \label{_conti_map_mero_Corollary_} 
Let $\phi:\; X \arrow Y$ be a continous map of complex analytic
varieties, which is holomorphic on an open dense subset 
$U\subset X$. Then $\phi$ is meromorphic.

\hfill

{\bf Proof:} Take a normalization $\tilde X \stackrel n \arrow X$.
Applying \ref{_conti_map_norma_holo_Corollary_} to the composition
$n\circ \phi:\; \tilde X \arrow Y$, we obtain that
$n\circ \phi$ is holomorphic. Thus, $\phi$ is meromorphic
as a composition of holomorphic
$n\circ \phi$ and a meromorphic map $n^{-1}$. \endproof

\hfill

Applying \ref{_conti_map_mero_Corollary_} to the situation
of \ref{_triana_subse_comple_ana_Theorem_}, we obtain the
following corollary.

\hfill

\corollary \label{_Psi_bimero_Corollary_}
In the situation of \ref{_triana_subse_isome_Theorem_},
the map \[ \Psi^{t_1}_{t_2}:\; X_{t_1} \arrow X_{t_2} \]
is bimeromorphic for each induced complex structures.

\endproof

\hfill

It remains to make a leap from ``bimeromorphic'' to ``holomorphic''.
This is done in two steps. We prove a number of
algebro-geometric statements about the behaviour of $\Psi^{t_1}_{t_2}$,
concluding with \ref{_pi_i_properties_Proposition_}.
In \ref{_alge_geo_suffi_pi_i_iso_Proposition_}, 
we show that these statements are strong enough to show that
the $\Psi^{t_1}_{t_2}$ is holomorphic. The premise
of \ref{_alge_geo_suffi_pi_i_iso_Proposition_} 
is purely algebro-geometric and its proof is independent
from the rest of this section.

We finish this Subsection with the following statement, which 
we use in \ref{_pi_i_properties_Proposition_}.

\hfill

\proposition \label{_Psi_iso_on_norma_Proposition_}
The map $\Psi^{t_1}_{t_2}:\; X_{t_1} \arrow X_{t_2}$ induces
an isomorphism of normalizations
\[ \tilde \Psi^{t_1}_{t_2}:\;\tilde X_{t_1} \arrow \tilde X_{t_2}. \]

{\bf Proof:} From the definition of normalization \cite{_Grauert_}, the
following lemma is evident.

\hfill

\lemma
Let $\Psi:\; X \arrow Y$ be a homeomorphism of complex
varieties. Assume that in an open dense subset $U\subset X$,
$\Psi$ is holomorphic. Then $\Psi$
induces an isomorphism $\tilde \Psi:\; \tilde X \arrow \tilde Y$ 
of normalizations, if the following statement holds.

\begin{description}
\item[(*)]  There exist a Stein covering $\{ U_i\}$ of $X$ such that
$\{ \Psi(U_i)\}$ is a Stein covering for $Y$.
\end{description}
\endproof

To prove \ref{_Psi_iso_on_norma_Proposition_}, it remains to show that
the property {\bf (*)} holds for the map 
$\Psi^{t_1}_{t_2}:\; X_{t_1} \arrow X_{t_2}$. 
For each point $x \in X_{t_1}$, it suffices to 
construct a Stein neighbourhood $U$ of $X_{t_1}$, such that 
$\Psi^{t_1}_{t_2}(U)$ is Stein. For a K\"ahler variety,
consider an open ball $B$ of radius $r$, taken with respect to the 
metric defined by geodesics. Then, for $r$ sufficiently small,
$B$ is Stein (\cite{_Greene_}). Since $\Psi^{t_1}_{t_2}$ is an 
isometry (\ref{_triana_subse_isome_Theorem_}), 
an image of an open ball of radius $r$ is again 
an open ball of radius $r$. This gives a 
system of Stein neighbourhoods satisfying
{\bf (*)}. \ref{_Psi_iso_on_norma_Proposition_}
is proven. \endproof


\subsection{Homeomorphisms of completely geodesic subvarieties 
induce isomorphisms of Zariski tangent spaces.}


\hfill

\proposition\label{_homeo_indu_iso_Zariski_Proposition_}
Let $M_1$, $M_2$ be a K\"ahler manifolds,
$X_1^{ns}\subset M_1$, $X_2^{ns}\subset M_2$ 
be completely geodesic complex submanifolds, not necessarily
closed, and $X_1$, $X_2$ be the closures of $X_1^{ns}$, $X_2^{ns}$ in
$M_1$, $M_2$. Assume that $X_1$, $X_2$ are complex analytic 
subvarieties of $M_1$, $M_2$. Let $\phi:\; X_1\arrow X_2$ be a morphism
of complex varieties, such that $\phi\restrict{X_1^{ns}}$ is an isometry.
Then $\phi$ induces an isomorphism of Zariski tangent spaces. 

\hfill

{\bf Proof:}
To prove \ref{_homeo_indu_iso_Zariski_Proposition_}, we interpret
the Zariski tangent space 
$T_x X_1$ in terms of the metric structure on $X_i^{ns}$, where $x$ 
is a point of $X_1$. Let $\gamma:\; [0,1] \arrow X_1$ be a path satisfying 
$\gamma\left([0,1] \backslash\{0.5\}\right) \subset X_1^{ns}$,
$\gamma(0.5) =x$. 
There is a natural topology on the total space $Tot(T X_1^{ns})$ ,
which comes from the embedding 
$Tot(TX_1) \stackrel i\hookrightarrow Tot(TM_1)$.
Since $X_1$ is completely geodesic, $i$ is an isometry.
This topology is compatible with the map 
$d\phi:\; TX_1 \arrow TX_2$, because $d\phi$ is also an isometry.
Assume that $\gamma$ is differentiable 
outside of $\{0.5\}$ and 
\[ \lim\limits_{t\to +0.5} \frac{d\gamma}{dt} = 
   \lim\limits_{t\to -0.5} \frac{d\gamma}{dt}.
\]
(the limits are taken in the metric completion of $TX_1$,
which might be considered as a subset of $Tot(TM_1)$).
The Zariski tangent space $T_x X_1$ can be identified with
equivalence classes of such paths as follows. Two 
paths $\gamma, \gamma'$ are equivalent if
\[ 
   \lim\limits_{t\to 0.5} 
   \frac{\rho(\gamma_1(t), \gamma_2(t))}{(t-0.5)^2} =0,
\]
where $\rho$ is the distance function in $M_1$. Since the
distance in $M_1$ coinsides with the distance in $X_1$, and
$\phi$ is isometry, this equivalence relation is compatible
with $\phi$. \endproof


\subsection{Algebro-geometric properties of the map $\Psi^{s_1}_{s_2}$.}


Let $M$ be a compact hyperk\"ahler manifold, $X\subset M$ a 
closed trianalytic subvariety. Fix a choice of induced complex structure.
Let $D(x)$ be the Douady space of $X$, and $s_1, s_2$ points
on $D(X)$, and $\gamma$ a real analytic path in $D(x)$
connecting $s_1$ and $s_2$. Let $X_1 = X_{s_1}, X_2 = X_{s_2}$
be the subvarieties of $M$ corresponding to $s_1$, $s_2$,
and $\Psi:\; X_1 \arrow X_2$ be the bimeromorphic map
of \ref{_Psi_bimero_Corollary_}. Let $X_1^{ns}$, $X_2^{ns}$ be the
non-singular part of $X_1$, and $\Psi^{ns}$ be the restriction of $\Psi$
to $X_1^{ns}$. Let $\Gamma^{ns}\subset M\times M$ be a graph of
$\Phi^{ns}$, and $\Gamma$ its closure. Clearly, $\Gamma$
is a complex analytic subvariety of $M\times M$.
Let $\pi_i:\; \Gamma \arrow X_i$ be the projections of
$\Gamma$ to $X_i$, which are obviously morphisms
of complex varieties. 

\hfill

\proposition\label{_pi_i_properties_Proposition_} 
The maps $\pi_i$, $i= 1,2$, have the following properties.

\begin{description}
\item[(i)] $\pi_i$ is finite.
\item[(ii)] $\pi_i$ is dominant and induces isomorphism
of normalizations
\item[(iii)] For every point $x\in \Gamma$, the differential
$d \pi_i :\; T_x \Gamma \arrow T_{\pi_i(x)} X_i$ is an
isomorphism.
\end{description}

{\bf Proof:} The statement 
(ii) follows from \ref{_Psi_iso_on_norma_Proposition_}
and (iii) from \ref{_homeo_indu_iso_Zariski_Proposition_}. 
To prove \ref{_pi_i_properties_Proposition_} (i), consider the
normalization map $n:\; \tilde \Gamma \arrow \Gamma$. On the
level of rings of functions, we have an embedding 
\[ 
   \calo(X_i) \hookrightarrow \calo(\Gamma) \hookrightarrow
   \calo(\tilde \Gamma).
\]
Since $\tilde \Gamma$ is a normalization of $X_i$ by (ii),
the ring $\calo(\tilde \Gamma)$ is finitely generated
as a $\calo(X_i)$-module. Since $\calo(X_i)$ is Noetherian
(\cite{_Grauert_}), this implies that $\calo(\Gamma)$
is also finitely generated as a $\calo(X_i)$-module.
This proves \ref{_pi_i_properties_Proposition_} (i).
\endproof

\hfill

To prove that the map $\Psi:\; X_1 \arrow X_2$, a.k.a.
$\Psi^{s_1}_{s_2}:\; X_{s_1} \arrow X_{s_2}$
is holomorphic, it suffices to show that $\pi_i:\; \Gamma \arrow X_i$
is an isomorphism. In the following Subsection, we show that
conditions (i) -- (iii) of \ref{_pi_i_properties_Proposition_} 
are {\it a priori} sufficient to establish that $\pi$ is 
an isomorphism. This will finish the proof of 
\ref{_triana_subse_comple_ana_Theorem_}.


\subsection{Finite dominant unramified morphisms of complex varieties.}


\hfill

\proposition \label{_alge_geo_suffi_pi_i_iso_Proposition_}
Let $\phi:\; X \arrow Y$ be a map of complex varieties satisfying
conditions (i)--(iii) of \ref{_pi_i_properties_Proposition_}.
Then $\phi$ is an isomorphism.

\hfill

{\bf Proof:} The map $\phi$ is one-to-one in general point (by (ii)). 
Thus, to prove that $\phi$ is an isomorphism it suffices to show
that $\phi$ is etale. On the other hand, $\phi$ is unramified by
(iii). By definition of etale morphisms,
to prove that $\phi$ is etale
it remains to show that $\phi$ is flat. To conclude the
proof of \ref{_alge_geo_suffi_pi_i_iso_Proposition_}, we
use the following lemma.

\hfill

\lemma \label{_unrami_domi_flat_Lemma_}
Let $\phi:\; X \arrow Y$ be a dominant morphism of complex varieties.
Assume that for every point $x\in X$, the map $\phi$ induces an
isomorphism $d\phi:\; T_x X \arrow T_{\phi(x)} Y$ of Zariski
tangent spaces. Then $\phi$ is flat.

\hfill

{\bf Proof:} Let $y=\phi (x)$. Conside the associated morphism
of local rings $\calo_y Y \stackrel {\phi_x} \hookrightarrow 
\calo_x X$. To prove
that $\phi$ is flat, it suffices to show that $\phi$ is an isomorphism.
Let ${\goth m}_x$, ${\goth m}_y$ be the maximal ideals
in $\calo_x X$, $\calo_y Y$. Then
${\goth m}_x/{\goth m}_x^2$ is generated by $\phi_x({\goth m}_y)$.
Thus, by Nakayama, $\phi_x({\goth m}_y)$ generate ${\goth m}_x$,
and we obtain ${\goth m}_x= \phi_x({\goth m}_y)\otimes \calo_x X$.
Consider $\calo_x X$ as $\calo_y Y$-module. Then 
$\calo_x X/{\goth m}_y \calo_x X$ is one-dimensional.
Applying Nakayama once more, we obtain that $\calo_x X$
is an $\calo_y Y$-module generated by $1\in\calo_x X$.
Thus, the map $\phi_x$ is surjective. Since $\phi$ is 
dominant, $\phi_x$ is also injective. Thus, $\calo_x X$
is a free $\calo_{\phi(x)} Y$-module for every $x\in X$.
By the local criterion of flatness, this implies that
$\phi$ is flat. \ref{_unrami_domi_flat_Lemma_} is proven.
This finishes the proof of \ref{_alge_geo_suffi_pi_i_iso_Proposition_}
and \ref{_triana_subse_comple_ana_Theorem_} \endproof


\section{Singular hyperk\"ahler varieties.}
\label{_singu_hype_Section_}


It is an intriguing question, what is the 
``correct''\footnote{``Correct'' in Platonic sense: some mathematicians
presume that the unique ``correct'' definition of each and every
significant mathematical object exists in itself and independently
of human perception.}
 definition of a singular hyperk\"ahler variety. We don't pretend to
answer this question. Instead, we give an {\it ad hoc} set of axioms
which describe some known examples (deformation spaces of stable bundles
and trianalytic subvarieties). It is likely that this {\it ad hoc}
definition is stronger than the ``correct'' one. A more elegant approach
was suggested by Deligne and Simpson 
(\cite{_Deligne:defi_}, \cite{_Simpson:hyperka-defi_}).

\hfill

\definition\label{_singu_hype_Definition_}
(\cite{_Verbitsky:Hyperholo_bundles_}, Definition 6.5)
Let $M$ be a hypercomplex variety (\ref{_hypercomplex_Definition_}).
The following data define a structure of {\bf hyperk\"ahler variety}
on $M$.

\begin{description}

\item[(i)] For every $x\in M$, we have an $\R$-linear 
symmetric positively defined
bilinear form $s_x:\; T_x M \times T_x M \arrow \R$
on the corresponding real Zariski tangent space.

\item[(ii)] For each triple of induced complex structures
$I$, $J$, $K$, such that $I\circ J = K$, we have a 
holomorphic differential 2-form $\Omega\in \Omega^2(M, I)$.

\item[(iii)] 
Fix a triple of induced complex structure 
$I$, $J$, $K$, such that $I\circ J = K$. Consider the
corresponding differential 2-form $\Omega$ of (ii).
Let $J:\; T_x M \arrow T_x M$ be an endomorphism of 
the real Zariski tangent spaces defined by $J$, and $Re\Omega\restrict x$
the real part of $\Omega$, considered as a bilinear form on $T_x M$.
Let $r_x$ be a bilinear form $r_x:\; T_x M \times T_x M \arrow \R$ 
defined by $r_x(a,b) = - Re\Omega\restrict x (a, J(b))$.
Then $r_x$ is equal to the form $s_x$ of (i). In particular,
$r_x$ is independent from the choice of $I$, $J$, $K$.

\end{description}

\noindent \remark \nopagebreak
\begin{description}
\item[(a)] It is clear how to define a morphism of hyperk\"ahler varieties.

\item[(b)]
For $M$ non-singular,  \ref{_singu_hype_Definition_} is
 equivalent to the usual
one (\ref{_hyperkahler_manifold_Definition_}). 
If $M$ is non-singular,
the form $s_x$ becomes the usual Riemann form, and 
$\Omega$ becomes the standard holomorphically symplectic form.

\item[(c)] It is easy to check the following.
Let $X$ be a hypercomplex subvariety of a hyperk\"ahler
variety $M$. Then, restricting the forms $s_x$ and $\Omega$
to $X$, we obtain a hyperk\"ahler structure on $X$. In particular,
trianalytic subvarieties of hyperk\"ahler manifolds are always
hyperk\"ahler, in the sense of \ref{_singu_hype_Definition_}.

\end{description}

\hfill

{\bf Caution:} Not everything which is seemingly hyperk\"ahler
satisfies the conditions of \ref{_singu_hype_Definition_}.
Take a quotient $M/G$ os a hyperk\"ahler manifold by an action 
of finite group $G$, acting in accordance with hyperk\"ahler
structure. Then $M/G$ is, generally speaking, {\it not} hyperk\"ahler
(in fact, $M/G$ is {\it never} hyperk\"ahler). For instance,
take a quotient of a 2-dimensional torus  $T$ by $G=\{\pm 1\}$
acting as an involution $t\arrow -t$. This is a beautiful and well
known example of a hyperk\"ahler automorphism; the quotient space
has 16 isolated singular points, which, if blown up, give
a K3 surface. For $x$ a singular point of $T/\{\pm 1\}$,
its Zariski tangent space has real dimension 6. 
On the other hand, for a hypercomplex variety,
there is a quaternion action in every Zariski tangent
space, and thus, the dimension real dimension
of Zariski tangent space must be divisible by 4. We obtain that
the space $T/\{\pm 1\}$ is not even hypercomplex. How this happens?
We take a twistor space $\Tw(T)$ of $T$ and take a quotient of 
$\Tw(T)$ by the natural action of $G=\{\pm 1\}$, which is 
holomorphic. The quotient is a complex variety fibered over $\C P^1$
For $\Tw(T)/\{\pm 1\}$ to be hypercomplex, this fibration must be trivial,
in a real analytic category. But the functor of forgetting the complex
structure does not commute with taking finite quotients! Thus, even
if $Tw(T)$ is (as a real analytic space) trivially fibered over
$\C P^1$, there is no way to push down this trivialization
to $\Tw(T)/\{\pm 1\}$. 

\hfill

The following theorem, proven in 
\cite{_Verbitsky:Hyperholo_bundles_} (Theorem 6.3), 
gives a convenient way to construct
examples of hyperk\"ahler varieties.

\hfill

\theorem \label{_hyperho_defo_hyperka_Theorem_}
Let $M$ be a compact hyperk\"ahler manifold, $I$ an induced
complex structure and $B$ a stable holomorphic bundle over $(M, I)$.
Let $D(B)$ be a deformation space of stable holomorphic structures on $B$.
Assume that $c_1(B)$, $c_2(B)$ are $SU(2)$-invariant, with respect
to the standard action of $SU(2)$ on $H^*(M)$. Then $D(B)$ has a
natural structure of a hyperk\"ahler variety. 

\nopagebreak
\endproof

\hfill

The following theorem is implicit in 
\cite{_Verbitsky:Hyperholo_bundles_}.

\hfill

\theorem \label{_hyperholo_functo_Theorem_}
Let $M$ be a compact hyperk\"ahler manifold, $I$ an induced
complex structure and $B_1, B_2, ..., B_n$ 
stable holomorphic bundles over $(M, I)$.
Let $D(B_i)$ be a deformation space of stable 
holomorphic structures on $B_i$.
Assume that $c_1(B_i)$, $c_2(B_i)$, $i=1, 2, ..., n$ 
are $SU(2)$-invariant, with respect
to the standard action of $SU(2)$ on $H^*(M)$. 
Let $\Pi$ be a natural tensor operation on the vector bundles,
such that, e. g.,

\[  B_1, ... B_n \arrow B_1 \otimes B_2 \otimes \Lambda^2 B_3 \otimes S^7
    B_4 \otimes ... \otimes B_n ^*.
\]
Assume that $\Pi(B_1, ... B_n)$ cannot be decomposed to a direct sum
of holomorphic bundles. Then $\Pi(B_1, ... B_n)$ is stable, 
and the associated map
\[ D(B_1) \times D(B_2) \times ..., \times D(B_n)\arrow 
   D\left(\Pi(B_1, ..., B_n)\right)
\]
(defined in a certain neighbourhood of 
$[B_1]\times [B_2]\times ..., \times [B_n] \in 
D(B_1) \times D(B_2) \times ..., \times D(B_n)$
is a morphism of hyperk\"ahler varieties.

\endproof

\hfill

\ref{_hyperholo_functo_Theorem_} gives a natural way to 
construct trianalytic subvarieites of hyperk\"ahler varieties. 

\hfill

The following theorem is almost trivial; the reader is advised 
to invent his or her own proof instead of reading ours
(which is by necessity sketchy).

\hfill

\theorem \label{_Doua_hyperka_Theorem_}
Let $M$ be a compact hyperk\"ahler manifold, $X\subset M$
a trianalytic subvariety and $D(X)$ its Douady 
space. Then $D(X)$ is a compact hyperk\"ahler variety.

\hfill

{\bf Proof:} Clearly, $D(X)$ is hypercomplex, in a natural way
(\ref{_Douady_hyperc_Remark_}). It remains to construct
the forms $s_x$ and $\Omega$. Let $Y$ be a trianalytic 
subvariety of $M$ which is a deformation of $X$,
and $U$ be a sufficiently small neighbourhood of $[Y]\in D(X)$. 
Let 
$\left(\pi:\; {\c X}_U \arrow U,\phi:\; {\c X}_U \arrow M\right)$ 
be the universal family of subvarieties
of $M$, attached to $U\subset D(X)$. The space ${\c X}_U$
is hypercomplex and the map  $\pi$
is compatible with the hypercomplex structure. 
Moreover, the map $\phi:\; {\c X}_U \arrow M$ is an
immersion, so ${\c X_U}$ is hyperk\"ahler (hyperk\"ahler structure
is obtained as a pullback from $M$).
\ref{_conne_in_triana_flat_Theorem_}
provides a natural trivialization of 
${\c X}_U \arrow U$. Thus, for each $y\in Y$, 
there exists a natural section $\sigma_y:\; U\hookrightarrow {\c X}_U$.
By \ref{_triana_subse_comple_ana_Theorem_}, 
this section is compatible with the hypercomplex structure.
Restricting the hyperk\"ahler structure
from ${\c X_U}$ to $\sigma_y(U)$, we obtain a hyperk\"ahler
structure on $U$. It is easy to check that this hyperk\"ahler
structure is independent from the choice of section
$\sigma_y$. Gluing the hyperk\"ahler structures from
different $U$, we obtain the proof of
\ref{_Doua_hyperka_Theorem_}. \endproof

\hfill

{\bf Acknowledgements:} I am indebted to D. Kaledin for enlightening
discussions; he also suggested some of the proofs. 
I am thankful to Tony Pantev for the frutiful discussions.
Mohan Ramachandran communicated me the reference
to \cite{_Greene_}. I am grateful to A. Beilinson and J. Bernstein 
for the help with algebraic geometry. My gratitude to
M. Finkelberg, D. Kazhdan and S.-T. Yau for their 
interest.

\hfill

\end{document}